\expandafter\chardef\csname pre amssym.def at\endcsname=\the\catcode`\@
\catcode`\@=11

\def\undefine#1{\let#1\undefined}
\def\newsymbol#1#2#3#4#5{\let\next@\relax
 \ifnum#2=\@ne\let\next@\msafam@\else
 \ifnum#2=\tw@\let\next@\msbfam@\fi\fi
 \mathchardef#1="#3\next@#4#5}
\def\mathhexbox@#1#2#3{\relax
 \ifmmode\mathpalette{}{\m@th\mathchar"#1#2#3}%
 \else\leavevmode\hbox{$\m@th\mathchar"#1#2#3$}\fi}
\def\hexnumber@#1{\ifcase#1 0\or 1\or 2\or 3\or 4\or 5\or 6\or 7\or 8\or
 9\or A\or B\or C\or D\or E\or F\fi}

\font\tenmsa=msam10
\font\sevenmsa=msam7
\font\fivemsa=msam5
\newfam\msafam
\textfont\msafam=\tenmsa
\scriptfont\msafam=\sevenmsa
\scriptscriptfont\msafam=\fivemsa
\edef\msafam@{\hexnumber@\msafam}
\mathchardef\dabar@"0\msafam@39
\def\dashrightarrow{\mathrel{\dabar@\dabar@\mathchar"0\msafam@4B}}
\def\dashleftarrow{\mathrel{\mathchar"0\msafam@4C\dabar@\dabar@}}

\def\ulcorner{\delimiter"4\msafam@70\msafam@70 }
\def\urcorner{\delimiter"5\msafam@71\msafam@71 }
\def\llcorner{\delimiter"4\msafam@78\msafam@78 }
\def\lrcorner{\delimiter"5\msafam@79\msafam@79 }
\def\yen{{\mathhexbox@\msafam@55 }}
\def\checkmark{{\mathhexbox@\msafam@58 }}
\def\circledR{{\mathhexbox@\msafam@72 }}
\def\maltese{{\mathhexbox@\msafam@7A }}

\font\tenmsb=msbm10
\font\sevenmsb=msbm7
\font\fivemsb=msbm5
\newfam\msbfam
\textfont\msbfam=\tenmsb
\scriptfont\msbfam=\sevenmsb
\scriptscriptfont\msbfam=\fivemsb
\edef\msbfam@{\hexnumber@\msbfam}

\catcode`\@=\csname pre amssym.def at\endcsname

\expandafter\ifx\csname pre amssym.tex at\endcsname\relax \else \endinput\fi
\expandafter\chardef\csname pre amssym.tex at\endcsname=\the\catcode`\@
\catcode`\@=11
\newsymbol\boxdot 1200
\newsymbol\boxplus 1201
\newsymbol\boxtimes 1202
\newsymbol\square 1003
\newsymbol\blacksquare 1004
\newsymbol\centerdot 1205
\newsymbol\lozenge 1006
\newsymbol\blacklozenge 1007
\newsymbol\circlearrowright 1308
\newsymbol\circlearrowleft 1309
\undefine\rightleftharpoons
\newsymbol\rightleftharpoons 130A
\newsymbol\leftrightharpoons 130B
\newsymbol\boxminus 120C
\newsymbol\Vdash 130D
\newsymbol\Vvdash 130E
\newsymbol\vDash 130F
\newsymbol\twoheadrightarrow 1310
\newsymbol\twoheadleftarrow 1311
\newsymbol\leftleftarrows 1312
\newsymbol\rightrightarrows 1313
\newsymbol\upuparrows 1314
\newsymbol\downdownarrows 1315
\newsymbol\upharpoonright 1316
 
\newsymbol\downharpoonright 1317
\newsymbol\upharpoonleft 1318
\newsymbol\downharpoonleft 1319
\newsymbol\rightarrowtail 131A
\newsymbol\leftarrowtail 131B
\newsymbol\leftrightarrows 131C
\newsymbol\rightleftarrows 131D
\newsymbol\Lsh 131E
\newsymbol\Rsh 131F
\newsymbol\rightsquigarrow 1320
\newsymbol\leftrightsquigarrow 1321
\newsymbol\looparrowleft 1322
\newsymbol\looparrowright 1323
\newsymbol\circeq 1324
\newsymbol\succsim 1325
\newsymbol\gtrsim 1326
\newsymbol\gtrapprox 1327
\newsymbol\multimap 1328
\newsymbol\therefore 1329
\newsymbol\because 132A
\newsymbol\doteqdot 132B
 
\newsymbol\triangleq 132C
\newsymbol\precsim 132D
\newsymbol\lesssim 132E
\newsymbol\lessapprox 132F
\newsymbol\eqslantless 1330
\newsymbol\eqslantgtr 1331
\newsymbol\curlyeqprec 1332
\newsymbol\curlyeqsucc 1333
\newsymbol\preccurlyeq 1334
\newsymbol\leqq 1335
\newsymbol\leqslant 1336
\newsymbol\lessgtr 1337
\newsymbol\backprime 1038
\newsymbol\risingdotseq 133A
\newsymbol\fallingdotseq 133B
\newsymbol\succcurlyeq 133C
\newsymbol\geqq 133D
\newsymbol\geqslant 133E
\newsymbol\gtrless 133F
\newsymbol\sqsubset 1340
\newsymbol\sqsupset 1341
\newsymbol\vartriangleright 1342
\newsymbol\vartriangleleft 1343
\newsymbol\trianglerighteq 1344
\newsymbol\trianglelefteq 1345
\newsymbol\bigstar 1046
\newsymbol\between 1347
\newsymbol\blacktriangledown 1048
\newsymbol\blacktriangleright 1349
\newsymbol\blacktriangleleft 134A
\newsymbol\vartriangle 134D
\newsymbol\blacktriangle 104E
\newsymbol\triangledown 104F
\newsymbol\eqcirc 1350
\newsymbol\lesseqgtr 1351
\newsymbol\gtreqless 1352
\newsymbol\lesseqqgtr 1353
\newsymbol\gtreqqless 1354
\newsymbol\Rrightarrow 1356
\newsymbol\Lleftarrow 1357
\newsymbol\veebar 1259
\newsymbol\barwedge 125A
\newsymbol\doublebarwedge 125B
\undefine\angle
\newsymbol\angle 105C
\newsymbol\measuredangle 105D
\newsymbol\sphericalangle 105E
\newsymbol\varpropto 135F
\newsymbol\smallsmile 1360
\newsymbol\smallfrown 1361
\newsymbol\Subset 1362
\newsymbol\Supset 1363
\newsymbol\Cup 1264
 
\newsymbol\Cap 1265
 
\newsymbol\curlywedge 1266
\newsymbol\curlyvee 1267
\newsymbol\leftthreetimes 1268
\newsymbol\rightthreetimes 1269
\newsymbol\subseteqq 136A
\newsymbol\supseteqq 136B
\newsymbol\bumpeq 136C
\newsymbol\Bumpeq 136D
\newsymbol\lll 136E
 
\newsymbol\ggg 136F
 
\newsymbol\circledS 1073
\newsymbol\pitchfork 1374
\newsymbol\dotplus 1275
\newsymbol\backsim 1376
\newsymbol\backsimeq 1377
\newsymbol\complement 107B
\newsymbol\intercal 127C
\newsymbol\circledcirc 127D
\newsymbol\circledast 127E
\newsymbol\circleddash 127F
\newsymbol\lvertneqq 2300
\newsymbol\gvertneqq 2301
\newsymbol\nleq 2302
\newsymbol\ngeq 2303
\newsymbol\nless 2304
\newsymbol\ngtr 2305
\newsymbol\nprec 2306
\newsymbol\nsucc 2307
\newsymbol\lneqq 2308
\newsymbol\gneqq 2309
\newsymbol\nleqslant 230A
\newsymbol\ngeqslant 230B
\newsymbol\lneq 230C
\newsymbol\gneq 230D
\newsymbol\npreceq 230E
\newsymbol\nsucceq 230F
\newsymbol\precnsim 2310
\newsymbol\succnsim 2311
\newsymbol\lnsim 2312
\newsymbol\gnsim 2313
\newsymbol\nleqq 2314
\newsymbol\ngeqq 2315
\newsymbol\precneqq 2316
\newsymbol\succneqq 2317
\newsymbol\precnapprox 2318
\newsymbol\succnapprox 2319
\newsymbol\lnapprox 231A
\newsymbol\gnapprox 231B
\newsymbol\nsim 231C
\newsymbol\ncong 231D
\newsymbol\diagup 231E
\newsymbol\diagdown 231F
\newsymbol\varsubsetneq 2320
\newsymbol\varsupsetneq 2321
\newsymbol\nsubseteqq 2322
\newsymbol\nsupseteqq 2323
\newsymbol\subsetneqq 2324
\newsymbol\supsetneqq 2325
\newsymbol\varsubsetneqq 2326
\newsymbol\varsupsetneqq 2327
\newsymbol\subsetneq 2328
\newsymbol\supsetneq 2329
\newsymbol\nsubseteq 232A
\newsymbol\nsupseteq 232B
\newsymbol\nparallel 232C
\newsymbol\nmid 232D
\newsymbol\nshortmid 232E
\newsymbol\nshortparallel 232F
\newsymbol\nvdash 2330
\newsymbol\nVdash 2331
\newsymbol\nvDash 2332
\newsymbol\nVDash 2333
\newsymbol\ntrianglerighteq 2334
\newsymbol\ntrianglelefteq 2335
\newsymbol\ntriangleleft 2336
\newsymbol\ntriangleright 2337
\newsymbol\nleftarrow 2338
\newsymbol\nrightarrow 2339
\newsymbol\nLeftarrow 233A
\newsymbol\nRightarrow 233B
\newsymbol\nLeftrightarrow 233C
\newsymbol\nleftrightarrow 233D
\newsymbol\divideontimes 223E
\newsymbol\varnothing 203F
\newsymbol\nexists 2040
\newsymbol\Finv 2060
\newsymbol\Game 2061
\newsymbol\mho 2066
\newsymbol\eth 2067
\newsymbol\eqsim 2368
\newsymbol\beth 2069
\newsymbol\gimel 206A
\newsymbol\daleth 206B
\newsymbol\lessdot 236C
\newsymbol\gtrdot 236D
\newsymbol\ltimes 226E
\newsymbol\rtimes 226F
\newsymbol\shortmid 2370
\newsymbol\shortparallel 2371
\newsymbol\smallsetminus 2272
\newsymbol\thicksim 2373
\newsymbol\thickapprox 2374
\newsymbol\approxeq 2375
\newsymbol\succapprox 2376
\newsymbol\precapprox 2377
\newsymbol\curvearrowleft 2378
\newsymbol\curvearrowright 2379
\newsymbol\digamma 207A
\newsymbol\varkappa 207B
\newsymbol\Bbbk 207C
\newsymbol\hslash 207D
\undefine\hbar
\newsymbol\hbar 207E
\newsymbol\backepsilon 237F
\catcode`\@=\csname pre amssym.tex at\endcsname

\magnification=1200
\hsize=468truept
\vsize=646truept
\voffset=-10pt
\parskip=4pt
\baselineskip=14truept
\count0=1

\dimen100=\hsize

\def\leftill#1#2#3#4{
\medskip
\line{$
\vcenter{
\hsize = #1truept \hrule\hbox{\vrule\hbox to  \hsize{\hss \vbox{\vskip#2truept
\hbox{{\copy100 \the\count105}: #3}\vskip2truept}\hss }
\vrule}\hrule}
\dimen110=\dimen100
\advance\dimen110 by -36truept
\advance\dimen110 by -#1truept
\hss \vcenter{\hsize = \dimen110
\medskip
\noindent { #4\par\medskip}}$}
\advance\count105 by 1
}
\def\rightill#1#2#3#4{
\medskip
\line{
\dimen110=\dimen100
\advance\dimen110 by -36truept
\advance\dimen110 by -#1truept
$\vcenter{\hsize = \dimen110
\medskip
\noindent { #4\par\medskip}}
\hss \vcenter{
\hsize = #1truept \hrule\hbox{\vrule\hbox to  \hsize{\hss \vbox{\vskip#2truept
\hbox{{\copy100 \the\count105}: #3}\vskip2truept}\hss }
\vrule}\hrule}
$}
\advance\count105 by 1
}
\def\midill#1#2#3{\medskip
\line{$\hss
\vcenter{
\hsize = #1truept \hrule\hbox{\vrule\hbox to  \hsize{\hss \vbox{\vskip#2truept
\hbox{{\copy100 \the\count105}: #3}\vskip2truept}\hss }
\vrule}\hrule}
\dimen110=\dimen100
\advance\dimen110 by -36truept
\advance\dimen110 by -#1truept
\hss $}
\advance\count105 by 1
}
\def\insectnum{\copy110\the\count120
\advance\count120 by 1
}

\font\ninerm=cmr9
\font\eightrm=cmr8

\font\tenrm=cmr10 at 10pt

\font\sc=cmcsc10

%\font\tenmsy=msym10
%\font\sevenmsy=msym7
%\font\fivemsy=msym5
%\newfam\msyfam
%\def\msy{\fam\msyfam\tenmsy}
%\textfont\msyfam=\tenmsy
%\scriptfont\msyfam=\sevenmsy
%\scriptscriptfont\msyfam=\fivemsy
\def\msb{\fam\msbfam\tenmsb}

\def\bbc{{\msb C}}

\def\bbi{{\msb I}}

\def\bbp{{\msb P}}
\def\bbq{{\msb Q}}
\def\bbr{{\msb R}}

\def\bbz{{\msb Z}}

\def\grD{\Delta}

\def\grG{\Gamma}

\def\grL{\Lambda}
\def\grO{\Omega}

\def\grS{\Sigma}

\def\grd{\delta}

\def\gri{\iota}

\def\grl{\lambda}

\def\gro{\omega}

\def\grt{\tau}

\def\grz{\zeta}

\def\la#1{\hbox to #1pc{\leftarrowfill}}
\def\ra#1{\hbox to #1pc{\rightarrowfill}}

\def\fract#1#2{\raise4pt\hbox{$ #1 \atop #2 $}}
\def\decdnar#1{\phantom{\hbox{$\scriptstyle{#1}$}}
\left\downarrow\vbox{\vskip15pt\hbox{$\scriptstyle{#1}$}}\right.}

\def\bowtie{\hbox to 1pt{\hss}\raise.66pt\hbox{$\scriptstyle{>}$}
\kern-4.9pt\triangleleft}
\def\hsmash{\triangleright\kern-4.4pt\raise.66pt\hbox{$\scriptstyle{<}$}}
\def\boxit#1{\vbox{\hrule\hbox{\vrule\kern3pt
\vbox{\kern3pt#1\kern3pt}\kern3pt\vrule}\hrule}}

\def\za{\vrule height6pt width4pt depth1pt}

\font\aa=eufm10

\def\Got#1{\hbox{\aa#1}}

\def\bfa{{\bf a}}

\def\bfw{{\bf w}}

\def\bfz{{\bf z}}

\def\calc{{\cal C}}
\def\calo{{\cal O}}

\def\cald{{\cal D}}

\def\calf{{\cal F}}

\def\calo{{\cal O}}

\def\cals{{\cal S}}

\def\calz{{\cal Z}}

\def\gF{{\Got F}}

\def\tU{\tilde{U}}

\font\svtnrm=cmr17

%\font\sbb=msym8
%\font\ninesl=cmsl9
\font\bsc=cmcsc10 at 10truept

\def\tU{\tilde{U}}

\def\mod{\hbox{mod}~}

\def\endo{\hbox{End}}
\def\lcm{\hbox{lcm}}
\def\codim{\hbox{codim}}

\def\exint{0}
\def\exsph{1}
\def\exsas{2}
\def\chfol{3}
\def\exfol{4}
\def\exproj{5}
\def\excon{6}
\def\exbr{7}

\centerline{\svtnrm Sasakian Geometry, Homotopy Spheres and}
\medskip
\centerline{\svtnrm Positive Ricci Curvature}
\bigskip
\centerline{\sc Charles P. Boyer~~ Krzysztof Galicki~~Michael Nakamaye}
\footnote{}{\ninerm During the preparation of this work the first two authors 
were partially supported by NSF grant DMS-9970904, and the third author by NSF grant 
DMS-0070190. 2000 Mathematics Subject Classification: 53C25, 57D60}
\bigskip

\bigskip\centerline{\vbox{\hsize = 5.85truein
\baselineskip = 12.5truept
\eightrm
\noindent {\bsc Abstract:}
We discuss the Sasakian geometry of odd dimensional homotopy spheres.  In
particular, we give a completely new proof of the existence of metrics
of positive Ricci curvature on exotic spheres that can be realized as the boundary of a 
parallelizable manifold. Furthermore, it is shown that on such homotopy spheres 
$\scriptstyle{\grS^{2n+1}}$ the moduli space of Sasakian structures has infinitely many 
positive 
components determined by inequivalent 
underlying contact structures. We also prove the existence of Sasakian metrics with 
positive Ricci curvature on each of the known $\scriptstyle{2^{2m}}$ distinct 
diffeomorphism types of homotopy real projective spaces 
$\scriptstyle{\bbr\bbp^{4m+1}}.$}} \tenrm

\bigskip
\bigskip
\bigskip
\bigskip
\bigskip
\baselineskip = 10 truept
\centerline{\bf \exint. Introduction}  
\bigskip
\bigskip

Milnor's discovery of exotic spheres [Mil1] presented Riemannian geometry with
a very natural question. What kind of special metrics or, more generally,
geometric structures can exist on exotic spheres? Perhaps the most intriguing
example of such a question concerns the existence of metrics with positive
sectional curvature. In 1974 Gromoll and Meyer [GM]
observed that one of the Milnor 7-spheres admits a metric of non-negative scalar 
curvature.Only recently Grove and Ziller [GZ] have observed that indeed all Milnor 
spheres, i.e., 7-spheres that are 3-sphere bundles over the 4-sphere admit metrics of
non-negative sectional curvature. Yet, it is not known whether any of these metrics
can be deformed to give a metric of positive sectional curvature.
However, it is known that some exotic spheres cannot admit such metrics.
This follows from a beautiful result of Hitchin, who
observed that  some of the exotic spheres, starting in dimension 9, do not even
admit metric of positive scalar curvature [Hi].

A somewhat more tractible problem concerns the existence of positive Ricci curvature
metrics on exotic spheres. Here, many examples have been known [Ch, Na, Her, Po].
In 1997, using surgery theory,
D. Wraith [Wr] proved the existence of Riemannian 
metrics with positive Ricci curvature on any exotic sphere that can be 
realized as the boundary of a parallelizable manifold. In particular, 
in dimension 7 all exotic spheres admit such metric. 

The question of the existence of other geometric structures on
exotic spheres has been equally important and intriguing. For example, 
it was realized mainly by the Japanese school [Abe1-2, AE, SH, Tak, Vai, YK]) as well as 
the French [LM] that Brieskorn manifolds naturally admit almost contact, contact, and 
Sasakian structures. The appoach closest to ours is that of Takahashi [Tak] who
constructed weighted Sasakian structures on Brieskorn manifolds. 

The purpose of this paper is to investigate the Sasakian geometry of homotopy 
spheres with the goal of producing Sasakian metrics with positive Ricci curvature on 
exotic spheres. The main ingredients of the proof are Brieskorn's description of homotopy 
spheres as links of isolated hypersurface singularities in $\bbc^{n+1},$ a positivity 
theorem in [BGN2] which is based on the orbifold or equivalently foliation version of the 
famous Calabi problem, and the algebraic geometry of Fano hypersurfaces in non 
well-formed weighted projective spaces. The link is the total space of 
an $S^1$ orbifold V-bundle 
over this Fano variety with a specific orbifold structure. From our point of view the 
differential invariants of the link are encoded in the orbifold structure, although the precise 
nature of this encoding still remains a mystery.

Recall that the classical Calabi Conjecture, later proved by Yau,
states that if $(M,\omega)$ is 
a smooth K\"ahler manifold then there exist a K\"ahler metric $\omega^\prime$ on $M$ 
whose
Ricci form $\rho^\prime$ has the property that $[\rho^\prime]=c_1(M)$. Moreover,
$[\omega^\prime]=[\omega]$, i.e., $\omega^\prime$ is in the same
cohomology class as $\omega$. In particular, in the Fano case
when $c_1(M)>0$ the $\omega^\prime$ has positive Ricci curvature. 
It has been often noted that the Calabi-Yau Conjecture
holds in the situation when $M$ is a K\"ahler orbifold (cf. [DK] or [Joy]). In the context of 
foliations  El Kacimi-Alaoui [ElK] actually gave a proof
of what one can call the ``transverse Yau theorem". In [BGN2] (see Theorem \exsas.10 
below) we showed that this theorem can be adapted to the situation,
where the foliation at hand is the characteristic foliation of  a Sasakian structure. 
In the present paper we combine all of these results to prove:

\noindent{\sc Theorem} A: \tensl For $n\geq 3$ let $\grS^{2n-1}$ be a homotopy 
sphere which can be realized as the boundary of a parallelizable manifold. 
Then $\grS^{2n-1}$ admits Sasakian metrics with positive Ricci curvature. 
\tenrm

We also show that there exists positive Sasakian structures belonging to infinitely 
many inequivalent contact structures. In particular, this implies that the moduli space of 
positive Sasakian structures on odd dimensional homotopy spheres has infinitely many 
components. Specifically, we prove the following:

\noindent{\sc Theorem} B: \tensl On each odd homotopy sphere 
$\grS^{2n-1}\in bP_{2n}$ there exists countably infinitely many deformation 
classes of positive Sasakian structures belonging to  non-isomorphic underlying  
contact  structures.  Hence, the moduli space of Sasakian structures on 
$\grS^{2n-1}$  has infinitely many positive components.  \tenrm

By a well-known theorem of Smale there are no exotic spheres in dimension 5; however,
there are 4 diffeomorphism types and 2 homeomorphism types of homotopy 
$\bbr\bbp^5$'s. More generally on the homotopy projective space $\bbr\bbp^{4m+1}$ 
there are at least $2^{2m}$ diffeomorphism types. Homotopy projective spaces have been 
studied in [AB,Bro,Gi1-2,LdM].  Our methods yield:
 
\noindent{\sc Theorem C}: \tensl On each of the known $2^{2m}$ diffeomorphism types of 
homotopy projective spaces $\bbr\bbp^{4m+1}$ there exists deformation types of 
positive Sasakian structures, and each deformation class contains Sasakian metrics of 
positive Ricci curvature. \tenrm

Finally, we employ a method of Savel'ev [Sav] to construct homotopy spheres from 
rational 
homology spheres. Combining this with previous work [BGN4] we give positive 
Sasakian 
structures on homotopy 9-spheres belonging to deformation classes of Sasakian 
structures that are inequivalent to the ones constructed previously using Brieskorn 
spheres.

\noindent{\sc Acknowledgments}: The authors would like to thank A. Buium, Y. Eliashberg, 
H. Hofer, and A. Vistoli for fruitful discussions on various aspects of this work, and Ian 
Hambleton for a clarification about homotopy projective spaces in the first version of this 
paper.

\vfill\eject
\bigskip
\baselineskip = 10 truept
\centerline{\bf \exsph. Exotic Spheres}  
\bigskip

Let us briefly recall a construction of exotic differential structures on odd
dimensional spheres. In 1956 Milnor stunned the mathematical world with the
construction of exotic differential structures on $S^7.$ Later the work of
Milnor and Kervaire [KeMi] (see also [Hz]) showed that associated with each
sphere $S^n$ with $n\geq 5$ there is an Abelian group $\Theta_n$ consisting of
equivalence classes of homotopy spheres $S^n$ that are equivalent under 
oriented h-cobordism. By Smale's h-cobordism theorem this implies equivalence 
under oriented diffeomorphism. The group operation on $\Theta$ is connected 
sum. $\Theta_n$ has a subgroup $bP_{n+1}$ consisting of equivalence classes of 
those homotopy $n$-spheres which bound parallelizable manifolds $V_{n+1}.$ 
Kervaire and Milnor [KeMi] proved that $bP_{2k+1}=0$  for $k\geq 1.$ Moreover, 
for $m\geq 2,$ $bP_{4m}$ is cyclic of order
$$|bP_{4m}|=2^{2m-2}(2^{2m-1}-1)~\hbox{numerator}~\bigl({4B_m\over m}\bigr),$$
where $B_m$ is the $m$-th Bernoulli number. For $bP_{4m+2}$ the situation is 
still not entirely understood. It entails computing the Kervaire invariant 
which is hard. It is known (see the recent review paper [La] and references 
therein) that  $bP_{4m+2}=0,$ or $\bbz_2,$ and is $\bbz_2$ if $4m+2\neq 2^i-2$ 
for any $i\geq 3.$ Furthermore, $bP_{4m+2}$ vanishes for $m=1,3,7,$ and $15.$

%In the case of $S^7$ we have
%$$\Theta_7=bP_8\approx \bbz_{28}.$$

Ten years after Milnor's landmark paper, Brieskorn [Br] showed how exotic
spheres can be obtained as links of isolated hypersurface singularities. In
particular if $g_m$ denotes the Milnor generator in $bP_{4m}$ (see [Hz]) then
Brieskorn showed that for each $k\in \bbz^+$ the link defined as the locus of 
points in $\bbc^{2m+1}$ that satisfies
$$\sum_{i=0}^{2m}|z_i^2|=1,\qquad z_0^{6k-1}+z_1^3+z_2^2+\cdots +z_{2m}^2=0,
\qquad k\geq 1 \leqno{\exsph.1}$$
represents an element in $bP_{4m}.$ We shall denote these homotopy
$(4m-1)$-spheres by $\grS_k^{4m-1}.$  The diffeomorphism type of $\grS\in 
bP_{4m}$ is determined [KeMi] by the signature $\grt$ of $V_{4m}$ which is 
necessarily divisible by 8. Two $\grS,\grS'\in bP_{4m}$ are diffeomorphic if 
and only if $\grt(V')\equiv \grt(V) (\mod 8|bP_{4m}|),$ and Brieskorn shows 
that $\grt(V_k)=(-1)^m8k,$ so that $\grS_k^{4m-1}$ and $\grS^{4m-1}_{k'}$ are 
diffeomorphic if and only if 
$$k'\equiv k (\mod |bP_{4m}|).$$    
The group operation in $bP_{4m}$ is connected sum and since the signature is 
additive under connected sum, the element $(-1)^m\grS_1^{4m-1}$ is a generator 
$g_m\in bP_{4m},$  and $(-1)^m|bP_{4m}|g_m$ is diffeomorphic to the
standard sphere $S^{4m-1}.$ 

More generally Brieskorn [Br,Di] considered the links $L(\bfa)$ defined by
$$\sum_{i=0}^{n}|z_i^2|=1,\qquad z_0^{a_0}+\cdots +z_{n}^{a_n}=0.
\leqno{\exsph.2}$$
To the vector $\bfa=(a_0,\cdots,a_n)\in \bbz_+^{n+1}$ one associates a graph 
$G(\bfa)$ whose $n+1$ vertices are labeled by $a_0,\cdots,a_n.$ Two vertices 
$a_i$ and $a_j$ are connected if and only if $\gcd(a_i,a_j)>1.$ Let $C_{ev}$ 
denote the connected component of $G(\bfa)$ determined by the even integers. 
Note that all even vertices belong to $C_{ev},$ but $C_{ev}$ may contain odd 
vertices as well. Then Brieskorn proved that $L(\bfa)$ is a $\bbz$-homology 
sphere if and only if either 
\item{1.} $G(\bfa)$ contains at least two isolated points, or 
\item{2.} $G(\bfa)$ contains one isolated point and $C_{ev}$ has an odd number 
of vertices and for any distinct $a_i,a_j\in C_{ev},~$ $\gcd(a_i,a_j)=2.$  

\noindent Hence, if $n\geq 3$ such an $L(\bfa)$ is homeomorphic to the sphere 
$S^{2n-1},$ and by Milnor's fibration theorem [Mil2] describes an element of 
$bP_{2n}.$

Clearly, the links in $bP_{4m}$ described by \exsph.1 are of type 1 above, 
whereas, examples of type 2 links in $bP_{4m+2}$ are given by 
$$\sum_{i=0}^{2m}|z_i^2|=1,\qquad z_0^p+z_1^2+z_2^2+\cdots +z_{2m+1}^2=0,
\qquad p\geq 1. \leqno{\exsph.3}$$
We shall denote these spheres by $\grS_p^{4m+1}.$ If $p\equiv \pm 1\mod 8$ 
then $\grS_p^{4m+1}$ is diffeomorphic to the standard sphere $S^{4m+1},$ 
whereas assuming that $4m+2\neq 2^i-2$ if $p\equiv \pm 3\mod 8$ then 
$\grS_p^{4m+1}$ is the exotic Kervaire sphere in $bP_{4m+2}.$ This follows 
from Levine's [Lev] computation of the Arf invariant [Br].

The polynomials given in \exsph.1 and \exsph.3 are by no means unique. Any 
other Brieskorn-Pham polynomial \exsph.2 satisfying either of the two 
conditions on $G(\bfa)$ above will do as well. For example, in leau of 
\exsph.1 we can consider 
$$\sum_{i=0}^{2m}|z_i^2|=1,\qquad z_0^{10k-1}+z_1^5+z_2^2+\cdots +z_{2m}^2=0,
\qquad k\geq 1 \leqno{\exsph.4}$$
In this case we find that $\grS^{4m-1}_k$ and $\grS^{4m-1}_{k'}$ are 
diffeomorphic if and only if $3k'\equiv 3k(\mod |bP_{4m}|).$

\bigskip
\baselineskip = 10 truept
\centerline{\bf \exsas. Sasakian Geometry on Homotopy Spheres}  
\bigskip

Recall [Bl,YK] that a Sasakian structure on a manifold $M$ of dimension $2n+1$ 
is a metric contact
structure $(\xi,\eta,\Phi,g)$ such that the Reeb vector field $\xi$ is a
Killing field and whose underlying almost CR structure is integrable. Briefly,
let $(M,\cald)$ be a contact manifold, and choose a 1-form $\eta$ so that
$\eta\wedge (d\eta)^n\neq 0$ and $\cald=\ker~\eta.$  The pair $(\cald,\gro)$,
where $\gro$ is the restriction of $d\eta$ to $\cald$ gives $\cald$ the
structure of a symplectic vector bundle. Choose an almost complex structure
$J$ on $\cald$ that is compatible with $\gro,$ that is $J$ is a smooth section
of the endomorphism bundle $\endo~\cald$ that satisfies 
$$ J^2= -\bbi, \qquad d\eta(JX,JY)=d\eta(X,Y), \qquad
d\eta(X,JX)>0\leqno{\exsas.1}$$ 
for any smooth sections $X,Y$ of $\cald.$ Notice
that $J$ defines a Riemannian metric $g_\cald$ on $\cald$
by setting $g_\cald(X,Y) =d\eta(X,JY).$ One easily checks that $g_\cald$
satisfies the compatibility condition $g_\cald(JX,JY)=g_\cald(X,Y).$ Now we
can extend $J$ to an endomorphism $\Phi$ on all of $TM$ by putting $\Phi =J$
on $\cald$ and $\Phi\xi=0.$ Likewise we can extend the metric $g_\cald$ on
$\cald$ to a Riemannian metric $g$ on $M$ by setting
$$g=g_\cald + \eta\otimes \eta. \leqno{\exsas.2}$$
The quadruple $(\xi,\eta,\Phi,g)$ is called a {\it metric contact structure}
on $M.$ If in addition $\xi$ is a Killing vector field and the almost complex
structure $J$ on $\cald$ is integrable the underlying almost contact
structure is said to be {\it normal} and $(\xi,\eta,\Phi,g)$ is called a 
{\it Sasakian structure}. The fiduciary examples of compact Sasakian
manifolds are the odd dimensional spheres $S^{2n+1}$ with the standard contact
structure and standard round metric $g.$

Actually as with K\"ahler structures there are many Sasakian structures on a
given Sasakian manifold. In fact there are many Sasakian structures 
which have $\xi$ as its characteristic vector field. If
$(\xi,\eta,\Phi,g)$ is a Sasakian structure on a smooth manifold (orbifold)
$M,$ we consider a deformation of this structure by adding to $\eta$ a 
continuous one parameter family of 1-forms $\grz_t$ that are basic with
respect to the characteristic foliation. We require that the 1-form
$\eta_t=\eta +\grz_t$  satisfy the conditions  
$$\eta_0=\eta, \qquad \grz_0=0,\qquad \eta_t\wedge (d\eta_t)^n\neq 0~~
\forall~~ t\in [0,1]. \leqno{\exsas.3}$$  
This last nondegeneracy condition implies that $\eta_t$ is a contact form on
$M$ for all $t\in [0,1]$ which by Gray's Stability Theorem 
belongs to the same underlying contact structure as $\eta.$ Moreover, since 
$\grz_t$ is basic $\xi$ is the Reeb (characteristic) vector field associated 
to $\eta_t$ for all $t.$ Now let us define
$$\eqalign{\Phi_t&=\Phi -\xi\otimes \grz_t\circ \Phi \cr
           g_t&=d\eta_t\circ(\Phi_t\otimes \hbox{id})+\eta_t\otimes \eta_t.} 
\leqno{\exsas.4}$$
In [BG3] it was proved that for all $t\in [0,1]$ and every basic 1-form 
$\grz_t$ such that $d\grz_t$ is of type $(1,1)$  and such that \exsas.3 holds
$(\xi,\eta_t,\Phi_t,g_t)$ defines a continuous 1-parameter family of Sasakian 
structures on $M$ belonging to the same underlying contact structure as 
$\eta.$ Given a Sasakian structure $\cals =(\xi,\eta,\Phi,g)$ on a manifold 
$M,$ we define $\gF(\xi)$ to be the family of all Sasakian structures obtained 
by the deformations above. 

Any Sasakian structure $\cals =(\xi,\eta,\Phi,g)$ has a 1-dimensional
foliation $\calf_\xi$ associated to it, defined by the flow of the Reeb vector 
field $\xi$ and called the {\it characteristic foliation}. Every Sasakian 
structure $\cals\in \gF(\xi)$ defines the same basic cohomology class 
$[d\eta]_B\in H^2_B(\calf)$ in the basic cohomology of the foliation 
$\calf(\xi),$ and conversely, any two homologous Sasakian structure $\cals 
=(\xi,\eta,\Phi,g)$ and $\cals' =(\xi,\eta',\Phi',g')$ lie in $\gF(\xi).$ 
Moreover, if $(M,\cals)$ is a compact Sasakian manifold then the cohomology 
class $[d\eta]_B$ is non-trivial, and cups to a non-trivial class in 
$H^{2n}_B(\calf_\xi).$ We are interested in the set of Sasakian structures 
which correspond to the same foliation $\calf_\xi.$ This set contains 
$\gF(\xi),$ but is slightly larger. We define $\gF(\calf_\xi)$ to be the set 
of all Sasakian structures whose characteristic foliation is $\calf_\xi.$ 
Clearly, we have $\gF(\xi)\subset \gF(\calf_\xi).$ For any Sasakian structure 
$\cals=(\xi,\eta,\Phi,g)$ there is the ``conjugate Sasakian structure'' 
defined by $\cals^c=(\xi^c,\eta^c,\Phi^c,g)=(-\xi,-\eta,-\Phi,g)\in 
\gF(\calf_\xi).$ Moreover, there is a type of homothety [YK] defined by
$$\xi'=a^{-1}\xi, \quad \eta'=a\eta, \quad \Phi'=\Phi, \quad 
g'=ag+(a^2-a)\eta\otimes \eta \leqno{\exsas.5}$$ 
for $a\in \bbr^+.$ So fixing $\cals$ we define
$$\gF^+(\calf_\xi)=\bigcup_{a\in \bbr^+}\gF(a^{-1}\xi), \leqno{\exsas.6}$$
and $\gF^-(\calf_\xi)$ to be the image of $\gF^+(\calf_\xi)$ under conjugation.
It is then easy to see that we have a decomposition 
$\gF(\calf_\xi)=\gF^+(\calf_\xi)\sqcup \gF^-(\calf_\xi).$
Notice that for $n$ even conjugation reverses orientation; whereas, for $n$ 
odd it preserves the orientation. This discussion shows that for each $\calf_\xi$             
the subset of homology classes represented by Sasakian structures forms a line 
minus the origin in $H^2_B(\calf_\xi).$ Since conjugation interchanges the positive and 
negative rays, we often restrict our considerations to the positive ray, and 
hence to $\gF^+(\calf_\xi)$. We shall refer to elements of 
$\gF^\pm(\calf_\xi)$ as {\it $a$-deformation classes} of Sasakian structures. 
We need

\noindent{\sc Definition} \exsas.7: \tensl Two Sasakian structures $\cals 
=(\xi,\eta,\Phi,g)$ and $\cals' =(\xi',\eta',\Phi',g')$ in $\gF(\calf_\xi)$ on 
a smooth manifold $M$ are said to be {\it $a$-homologous} if there is an $a\in 
\bbr^+$ such that $\xi'=a^{-1}\xi$ and $[d\eta']_B=a[d\eta]_B.$   \tenrm

The $a$-homology classes form a set of two elements that can be identified 
with positive and negative rays in $H^2_B(\calf_\xi).$ So every Sasakian 
structure in $\gF(\calf_\xi)$ is $a$-homologous to $\cals$ or its conjugate 
$\cals^c.$ 

Other important invariants of $\gF(\xi)$ are the basic Chern classes 
$c_i(\calf_\xi)$ of the symplectic vector bundle $(\cald,d\eta)$ as elements 
of the basic cohomology ring $H^*_B(\calf_\xi).$ In particular we are 
interested in the basic first Chern class $c_1(\calf_\xi)\in 
H^2_B(\calf_\xi).$ Recall [BGN2] that a Sasakian structure 
$\cals=(\xi,\eta,\Phi,g)$ is {\it positive (negative)} 
if its basic first Chern class $c_1(\calf_\xi)\in H^2(\calf_\xi)$ can be 
represented by a positive (negative) definite $(1,1)$-form. 
On rational homology spheres there are are only two possibilities, viz.

\noindent{\sc Proposition} \exsas.8: \tensl Let $\cals=(\xi,\eta,\Phi,g)$ be a 
 Sasakian structure on a rational homology sphere $M^{2n-1}.$ Then either 
$c_1(\calf_\xi)>0$ or $c_1(\calf_\xi)<0.$ \tenrm

\noindent{\sc Proof}: For rational homology spheres $M$ the basic long exact 
sequence
$$\cdots\ra{2.5}H_B^p(\calf_\xi)\ra{2.5}H^p(M,\bbr)\fract{j_p}{\ra{2.5}}
H_B^{p-1}(\calf_\xi) \fract{\grd}{\ra{2.5}} H^{p+1}_B(\calf_\xi)\ra{2.5}\cdots 
$$
implies the isomorphisms
$$H^{p-1}_B(\calf_\xi)\approx H^{p+1}_B(\calf_\xi)$$
for $p=1,\cdots,2n-1,$ and since $H^1_B(\calf_\xi)=0$ and $H^2_B(\calf_\xi)$ 
is generated by $[d\eta]_B,$ this implies the ring isomorphism
$$H^*_B(\calf_\xi)\approx \bbr\bigr[[d\eta]_B\bigl]/([d\eta]_B)^n. 
\leqno{\exsas.9}$$
Thus, $\calc(\calf_\xi)=\bbr^+\cup \bbr^-,$ and $c_1(\calf_\xi)=a[d\eta]_B$ 
for some $a\in \bbr.$ Since $d\eta$ is the transverse K\"ahler form it is 
positive definite, and  $a$ cannot vanish since this would imply that the 
basic geometric genus [BGN2] $p_g(\calf_\xi)\neq 0$ implying 
$b_{n-1}(M)\neq 0.$ Thus, either $a>0$ or $a<0$ implying either 
$c_1(\calf_\xi)>0$ or $c_1(\calf_\xi)<0.$  \hfill\za

In the next two sections we shall discuss how to construct positive Sasakian 
structures on homotopy spheres, and prove Theorem A of the Introduction.
Recall that a positive Sasakian structure [BGN2] is a Sasakian structure whose basic first 
Chern class can be represented by a positive definite basic $(1,1)$-form. Once positive 
Sasakian structures are obtained Sasakian structures with positive Ricci curvature follow 
from the following theorem of the authors.

\noindent{\sc Theorem} \exsas.10 [BGN2]: \tensl Let $\cals=(\xi,\eta,\Phi,g)$ 
be a positive Sasakian structure on a compact manifold $M$ of dimension 
$2n+1.$ Then $M$ admits a Sasakian structure $\cals'=(\xi',\eta',\Phi',g')$  
with positive Ricci curvature $a$-homologous to $\cals$ for some $a>0.$
\tenrm

\vfill\eject
\bigskip
\baselineskip = 10 truept
\centerline{\bf \chfol. The Characteristic Foliation and K\"ahler Orbifolds}  
\bigskip

This foliation encodes important
invariants of the Sasakian structure. For example, $(M,\xi,\eta,\Phi,g)$ is
said to be {\it quasi-regular} if all the leaves of $\calf_\xi$ are compact.
In this case the leaves are all circles, but in general they can have
nontrivial holonomy. If the leaf holonomy group is trivial for every leaf, the
Sasakian manifold $(M,\xi,\eta,\Phi,g)$ is called {\it regular}. It is well-known
that if $(M,\xi,\eta,\Phi,g)$ is compact and quasi-regular, then the
space of leaves $\calz$ is a compact K\"ahler orbifold, and a smooth manifold
in the regular case, whose K\"ahler class $[\gro]$ is represented by an
integral class in $H^2_{orb}(\calz,\bbz)$ and a rational class in 
$H^2(\calz,\bbq).$ The leaf holonomy groups of $\calf_\xi$
are the local uniformizing groups of the orbifold $\calz$ and these are
invariants of the Sasakian structure on $M.$ A somewhat courser invariant is
the {\it order} $\upsilon(M)$ of a compact Sasakian manifold defined to be the
least common multiple of the orders of the leaf holonomy groups. So
$(M,\xi,\eta,\Phi,g)$ is regular if and only if $\upsilon_M =1.$ 

For a quasi-regular Sasakian structure $\cals$ on a compact manifold $M$ the 
space of leaves $\calz$ is also a projective algebraic variety, and it is 
important to distinguish between $\calz$ as an orbifold and $\calz$ as an 
algebraic variety. In general the singular sets are different. In the case of 
the Brieskorn-Kervaire exotic spheres considered in Proposition \exfol.3 
below, we shall see that $\calz$ as an algebraic variety is smooth, in fact a 
projective space, whereas, its orbifold singular set is actually a divisor. For 
hypersurfaces (or more generally complete intersections) in weighted 
projective spaces, the two singular sets coincide precisely when the 
hypersurfaces are {\it well-formed} (cf. [Fl]); however, the orbifolds $\calz$ 
associated to the Brieskorn spheres are never well-formed. In the non well-formed
case certain pathologies can occur. For, example if $\calz\subset \bbp(w_0,\cdots, 
w_n)$ is not well-formed, there can be isomorphic coherent sheaves $\calo_\calz(m)$ 
and $\calo_\calz(n)$ with $n\neq m.$

In our previous work [BG1,BG2,BGN1,BGN2,BGN3] we have always worked with the 
well-formed case, and there was little need to distinguish between the 
orbifold context and the algebraic geometric context. In the present paper 
this is no longer the case. In algebraic geometry one uses a sort of poetic 
license not to distinguish between vector bundles and locally free sheaves. 
There is a similar situation in the category of orbifolds though much less 
well-known. First, the notion of  a vector V-bundle was introduced by Satake 
[Sat] and essentially consists of vector bundles defined on the local 
uniformizing covers $(\tU_i,\grG_i,\phi_i)$ that patch together in a nice way,  together 
with group homomorphisms $\grG_i\ra{1.3} GL(n,\bbc)$ that satisfy certain compatibility 
conditions (cf. [BG1] for details). On the other hand the notion of a locally V-free sheaf 
on a normal projective variety with only quotient singularities has been defined 
([Bl], see also [Ka] where it is referred to as a $\bbq$-vector bundle), and 
proves to be convenient in our context. A coherent sheaf $\gF$ of 
$\calo_\calz$-modules is {\it locally V-free} if on each uniformizing 
neighborhood $\tU_i$ there is a free sheaf $\tilde{\gF}$ together with an 
action of $\grG_i$ on $\tilde{\gF}$ such that 
$\gF=((\phi_i)_*\tilde{\gF})^{\grG_i}$ where $(\cdot)^\grG$ denotes the 
$\grG$-invariant subsheaf. Notice that a locally V-free sheaf is not 
necessarily locally free; however, a locally free sheaf is locally V-free. Now 
it is straightforward to see that the usual correspondence between vector 
V-bundles and locally V-free sheaves still holds. We shall also refer to a locally V-free 
sheaf of rank one as a {\it V-invertible sheaf}. V-invertible sheaves 
correspond to holomorphic line V-bundles. In our previous work [BG1,BG2] we 
introduced the group $\hbox{Pic}^{orb}(\calz)$ of holomorphic line V-bundles 
on $\calz.$ We also refer to $\hbox{Pic}^{orb}(\calz)$ as the group of 
V-invertible sheaves on $\calz.$ The {\it dualizing sheaf} $\gro_\calz$  
on $\calz,$ defined by
$$\gro_\calz = \gri_*(\grO^n\calz_{reg}),\leqno{\chfol.1}$$
where $\grO^p\calz_{reg}$ denotes the sheaf of differential p-forms on the 
algebro-geometric regular set $\calz_{reg},$ is of particular importance. For a given 
orbifold structure on $\calz$, $\gro_\calz$ is V-invertible with respect to that structure, 
and corresponds to the canonical V-bundle $K_\calz$ as defined by Baily [Ba].
However, $\calz$ may admit several orbifold structures, and $\gro_\calz$ can take 
different forms as V-invertible sheaves. 

Consider the affine space $\bbc^{n+1}$ together with a weighted
$\bbc^*$-action given by 
$$(z_0,\ldots,z_n) \mapsto 
(\grl^{w_0}z_0,\ldots,\grl^{w_n}z_n),\leqno{\chfol.2}$$ 
where the {\it weights} $w_j$ are
positive integers. It is convenient to view the weights as the components of a
vector $\bfw\in (\bbz^+)^{n+1},$ and we shall assume that they are ordered
$w_0\leq w_1\leq \cdots\leq w_n$ and that $\gcd(w_0,\ldots,w_n)=1.$ Let $f$
be a weighted homogeneous polynomial, that is $f\in \bbc[z_0,\ldots,z_n]$ and
satisfies 
$$f(\grl^{w_0}z_0,\ldots,\grl^{w_n}z_n)=\grl^df(z_0,\ldots,z_n),\leqno{\chfol.3}$$
where $d\in \bbz^+$ is the degree of $f.$ We are interested in the {\it
weighted affine cone} $C_f$ defined by
the equation $f(z_0,\ldots,z_n)=0.$ We shall assume that the origin in
$\bbc^{n+1}$ is an isolated singularity, in fact the only singularity, of
$f.$ Then the link $L_f$ defined by 
$$L_f= C_f\cap S^{2n+1},$$
where 
$$S^{2n+1}=\{(z_0,\ldots,z_n)\in \bbc^{n+1}|\sum_{j=0}^n|z_j|^2=1\}$$
is the unit sphere in $\bbc^{n+1},$ is a smooth manifold of dimension $2n-1.$ 
Furthermore, it is well-known [Mil2] that the link $L_f$ is $(n-2)$-connected.

On $S^{2n+1}$ there is a well-known [YK] ``weighted'' Sasakian structure  
$(\xi_\bfw,\eta_\bfw,\Phi_\bfw,g_\bfw)$ which in the standard coordinates
$\{z_j=x_j+iy_j\}_{j=0}^n$ on $\bbc^{n+1}=\bbr^{2n+2}$ is determined by
$$\eta_\bfw = {\sum_{i=0}^n(x_idy_i-y_idx_i)\over\sum_{i=0}^n
w_i(x_i^2+y_i^2)}, \qquad \xi_\bfw
=\sum_{i=0}^nw_i(x_i\partial_{y_i}-y_i\partial_{x_i}),$$
and the standard Sasakian structure $(\xi,\eta,\Phi,g)$ on $S^{2n+1}.$
The embedding $L_f\hookrightarrow S^{2n+1}$ induces a Sasakian structure on
$L_f$ [BG3]. 

Given a sequence $\bfw =(w_0,\ldots,w_n)$ of ordered positive integers one can
form the graded polynomial ring $S(\bfw)=\bbc[z_0,\ldots,z_n]$, where $z_i$ has
or {\it weight} $w_i.$ The weighted projective space [Dol, Fle]
$\bbp(\bfw)=\bbp(w_0,\ldots,w_n)$ is defined to be the scheme
$\hbox{Proj}(S(\bfw)).$  It is the quotient space 
$(\bbc^{n+1}-\{0\})/\bbc^*(\bfw)$, where $\bbc^*(\bfw)$ is the weighted action
defined in \chfol.2, or equivalently, $\bbp(\bfw)$ is the quotient of
the weighted Sasakian sphere
$S_\bfw^{2n+1}=(S^{2n+1},\xi_\bfw,\eta_\bfw,\Phi_\bfw,g_\bfw)$ by the
weighted circle action $S^1(\bfw)$ generated by $\xi_\bfw.$ As such
$\bbp(\bfw)$ is also a compact complex orbifold with an induced K\"ahler
structure. Furthermore, the hypersurface $\calz_f\subset \bbp(\bfw)$ defined as the 
zero 
set of the section $f\in \calo_{\bbp(\bfw)}(d)$ is the $S^1(\bfw)$ quotient of the link 
$L_f.$ 
The assumption that the origin of $\bbc^{n+1}$ is the only singularity of $f$ guarantees 
that the quotient space $\calz_f$ is an orbifold. In this case we also say the $\calz_f$ 
is {\it 
quasi-smooth}.  We have from [BG3]

\noindent{\sc Theorem} \chfol.4: \tensl The quadruple
$(\xi_\bfw,\eta_\bfw,\Phi_\bfw,g_\bfw)$ gives $L_f$ a quasi-regular Sasakian
structure $\cals_\bfw$ such that there is a commutative diagram
$$\matrix{L_f &\ra{2.5}& S^{2n+1}_\bfw&\cr
  \decdnar{\pi}&&\decdnar{} &\cr
   \calz_f &\ra{2.5} &\bbp(\bfw),&\cr}$$
where the horizontal arrows are Sasakian and K\"ahlerian embeddings,
respectively, and the vertical arrows are principal $S^1$ V-bundles and
orbifold Riemannian submersions.   Moreover, the orbifold first Chern
class in $c_1(\calz_f)\in H^2_{orb}(\calz_f,\bbz)$ is related to the basic 
first Chern class $c_1(\calf_{\xi_\bfw})$ of $\cals_\bfw$ by 
$c_1(\calf_{\xi_\bfw})=\pi^*c_1(\calz_f).$  \tenrm 

Hence, to prove Theorem A we need to compute the orbifold first Chern class 
$c_1(\calz_f).$ This requires determining the canonical V-bundle $K_{\calz_f}$ 
by an adjunction formula. The following is well-known [Dol,BR]:

\noindent{\sc Lemma} \chfol.5: [Dol] \tensl On $\bbp(\bfw)$ the dualizing sheaf 
$\gro_{\bbp(\bfw)}$ is isomorphic to $\calo_{\bbp(\bfw)}(-|\bfw|).$ \tenrm

\noindent{\sc Warning}: This does not imply that the anti-dualizing sheaf 
$\gro^{-1}_{\bbp(\bfw)}$ is isomorphic to $\calo_{\bbp(\bfw)}(|\bfw|).$ We note that this 
isomorphism does hold if the weighted projective space $\bbp(\bfw)$ is well-formed, 
i.e. 
the gcd of any n components of $\bfw=(w_0,\cdots,w_n)$ is one (cf. [BG3]), but here we 
are interested in the non well-formed case.

Next we are interested in an adjunction formula. This is known to hold when the 
hypersurface is well-formed as defined by Fletcher [Fl]; that is a hypersurface 
$\calz_f\subset \bbp(\bfw)$ is {\it well-formed} if $\bbp(\bfw)$ is well-formed and 
$\calz_f$ contains no codimension 2 singular stratum of $\bbp(\bfw).$  However, the 
condition that guarantees the adjunction formula is that the hypersurface $\calz_f$ 
contain no codimension 2 singular stratum of $\bbp(\bfw),$ and no further condition on 
$\bbp(\bfw)$ is needed. Indeed, in the sequel we shall make use of the adjunction 
formula for hypersurfaces in weighted projective spaces that are not well-formed.

\noindent{\sc Lemma} \chfol.6 [Adjunction formula]: \tensl Let $\calz_f$ be a 
quasi-smooth
hypersurface in $\bbp(\bfw)$ defined by a weighted homogeneous polynomial $f$ 
of degree $d.$ Suppose that $\codim(\calz_f\cap \bbp^{sing}_\bfw,\calz_f)\geq 2.$  
Then the adjunction formula  
$$\gro_{\calz_f}\approx \calo_{\calz_f}(d-|\bfw|)$$
holds. 
\tenrm

\noindent{\sc Proof}: By Proposition 5.75 of 
[KM], $\omega_X \simeq \calo_X(K_X)$
for any normal variety $X$.  
In particular, for the hypersurface $\calz_f$ in $\bbp(\bfw)$ we 
have
$$
\omega_{\calz_f} \approx \calo_{\calz_f}(K_{\calz_f}).
$$

In general, $\calz_f$ will not be a Cartier divisor and so we cannot
apply the adjunction formula directly.
On the other hand, if $U \subset \bbp(\bfw)$ denotes the smooth locus and
$U_f = \calz_f \cap U$ then the adjunction formula says that
$$
K_{U_f} = K_U + U_f|U_f = \calo_{U_f}(d-|\bfw|).
$$
Consequently
$$
\omega_{U_f} \approx \calo_{U_f}(d-|\bfw|)
$$
The lemma then follows from Proposition 0-1-10 of [KMM] since $\omega_{Z_f}$
and $\calo_{\calz_f}(d-|\bfw|)$ are reflexive rank one sheaves which agree off 
a set of codimension at least 2. \hfill\za

\noindent{\sc Remark} \chfol.7: In several of our applications the hypersurface $\calz_f$ is 
a Cartier divisor on $\bbp(\bfw)$, and in this case one can directly apply the adjunction 
formula in [KM].

The dualizing sheaf $\gro_{\calz_f}$ is a V-invertible sheaf with respect to the orbifold 
structure on $\calz_f,$ and as such corresponds to the canonical V-bundle 
$K_{\calz_f},$ 
or equivalently to the canonical  V-divisor of Baily [Ba], i.e.
$$\gro_{\calz_f}=\calo(K_{\calz_f}).\leqno{\chfol.8}$$
Now $K_{\calz_f}$ is an element of $\hbox{Pic}^{orb}(\calz_f)$ and as such has an 
inverse, namely, the anti-canonical V-bundle 
$K_{\calz_f}^{-1}.$ 
We recall [BG1] that a compact complex orbifold $\calz$ is {\it Fano} if its 
anti-canonical 
V-bundle $K_{\calz_f}^{-1}$ is ample, and in this case the {\it orbifold Fano index} of 
$\calz$ is defined to be the largest positive integer $I$
such that ${1\over I}K^{-1}_\calz$ is an element of $\hbox{Pic}^{orb}(\calz).$ Now the 
group $\hbox{Pic}^{orb}(\calz)$ is isomorphic [BG4] to the divisor class group 
$\hbox{Cl}(\calz).$ Since $\hbox{Cl}(\calz)$ is an algebraic invariant, the orbifold Fano 
index is also a purely algebraic invariant.   Furthermore, on a weighted projective space 
$\bbp(\bfw)$ we have that $\hbox{Cl}(\bbp(\bfw))\approx \bbz$ and $\hbox{Cl}(\bbp(\bfw))$ 
is generated by $\calo_{\bbp(\bfw)}(a_\bfw)$ [BR] where $a_\bfw=\lcm(d_0,\cdots,d_n)$ 
and $d_i=\gcd(w_0,\cdots,\hat{w_i},\cdots,w_n).$ We have

\noindent{\sc Lemma} \chfol.9: \tensl The following hold:
\item{a)}There is an isomorphism $\hbox{Pic}^{orb}(\bbp(\bfw)) \approx \bbz$ and as a 
V-bundle (or equivalently V-invertible sheaf) $\calo_{\bbp(\bfw)}(a_\bfw)$ is the 
generator of $\hbox{Pic}^{orb}(\bbp(\bfw))$ corresponding to $1\in \bbz.$ 
\item{b)} The (orbifold) Fano index $I$ is independent of the orbifold structure and on 
$\bbp(\bfw)$ equals $|\bar{\bfw}|$ where $\bar{\bfw}$ is the normalization of $\bfw.$
\item{c)} For any integers $m,n$ there is a V-isomorphism of V-invertible sheaves 
$$\calo_{\bbp(\bfw)}(na_\bfw)\otimes \calo_{\bbp(\bfw)}(ma_\bfw)\ra{1.8} 
\calo_{\bbp(\bfw)}((n+m)a_\bfw)$$ 
describing the group structure of $\hbox{Pic}^{orb}(\bbp(\bfw)).$ 
\item{d)} $\hbox{Pic}(\bbp(\bfw))$ is the subgroup of $\hbox{Pic}^{orb}(\bbp(\bfw))$ 
generated by $\calo_{\bbp(\bfw)}(\upsilon_\bfw)$ where $\upsilon_\bfw=\lcm(w_0,\cdots, 
w_n)$ is the order of $\bbp(\bfw)$ as an orbifold. \tenrm

\noindent{\sc Remark}. It is clear from the definitions that $a_\bfw$ divides 
$\upsilon_\bfw.$

Lemma \chfol.9 follows from the discussion above and the next lemma which is an 
orbifold version of a result of Delorme [Del].

\noindent{\sc Lemma} \chfol.10: \tensl For any $j\in \bbz$ there is a unique $J\in \bbz$ and 
a V-isomorphism of V-invertible sheaves
$$\calo_{\bbp(\bfw)}(j)\approx \calo_{\bbp(\bfw)}(Ja_\bfw).$$ 
\tenrm

\noindent{\sc Proof}: As an isomorphism of reflexive sheaves this is a result of Delorme 
[Del]. This isomorphism can be 
described as follows (cf. [BR]):  first notice that 
$\calo_{\bbp(\bfw)}(j)=(S(\bfw)(j))\widetilde{}$ where the graded $S(\bfw)$-module 
$S(\bfw)(j)$ is defined by $S(\bfw)(j)_l=S(\bfw)_{j+l}.$ Now there is an equality of 
schemes 
$\hbox{Proj}(S(\bfw)) = \hbox{Proj}(S'(\bfw))$ where $S'(\bfw)$ is the subring of 
$S(\bfw)$ 
defined by
$$S'(\bfw)=\bigoplus_{k\in \bbz}S(\bfw)_{ka_\bfw}.$$
The isomorphism in \chfol.10 is then induced by the equality of graded 
$S'(\bfw)$-modules
$$\bigoplus_{k\in \bbz}S(\bfw)(j)_{ka_\bfw} = \bigoplus_{k\in \bbz}z_0^{b_0(j)}\cdots 
z_n^{b_n(j)}\bigl(S(\bfw)(j-\sum_ib_i(j)w_i)\bigr)_{ka_\bfw}$$
where $b_i(j)$ are the unique integers $0\leq b_i(j)<d_i$ satisfying $j=b_i(j)w_i+c_id_i$ 
for some $c_i\in \bbz.$ So the integer $J$ is $j-\sum_ib_i(j)w_i.$

Let $\{U,\phi,\grG\}$ be a local 
uniformizing system. But since the isomorphic sheaves $\calo_{\bbp(\bfw)}(j)$ and 
$\calo_{\bbp(\bfw)}(Ja_\bfw)$ are locally V-free, they pull-back to isomorphic sheaves 
\break $\phi^*\calo_{\bbp(\bfw)}(j)^{\vee\vee}$ and 
$\phi^*\calo_{\bbp(\bfw)}(Ja_\bfw)^{\vee\vee}$ on the local cover $U.$  \hfill\za 

Our next task is to compute the orbifold Fano index.  It appears to be difficult to get a 
general formula, so we content ourselves with less. Actually all that is needed for this 
paper is the positivity of the index, but with a little extra effort we can do better. All of the 
weighted hypersurfaces considered in this paper are branched covers of the form 
$$g=z_0^p+f(z_1,\cdots,z_n)=0,\leqno{\chfol.11}$$
where $f$ is a weighted homogeneous polynomial of degree $d$ and weights 
$\bfw=(w_1,\cdots,w_n)$ that are reduced, but not necessarily well-formed. We also 
assume that $f$ is quasi-smooth and that $\gcd(d,p)=1.$ Under these conditions we have

\noindent{\sc Lemma} \chfol.12: \tensl Let $\calz_f$ be a quasi-smooth 
hypersurface in a weighted projective space $\bbp(\bfw).$ 
Suppose also that $\calz_f$ is cut out by the zero locus of the weighted 
homogeneous polynomial $f(z_1,\cdots, z_n)$ of degree $d.$ Let $p$ be a positive 
integer such that $\gcd(p,d)=1.$ Then the hypersurface $\calz_g$ of degree $d_g=dp$ cut 
out by equation \chfol.11 is quasi-smooth and Fano with orbifold Fano index $|\bar{\bfw}|$ 
where $\bar{\bfw}$ is the normalization of $\bfw,$ and as an algebraic variety 
$\calz_g$ is isomorphic to the weighted projective space $\bbp(\bar{\bfw}).$
\tenrm

\noindent{\sc Proof}: It is easy to see that $\calz_g$ is quasi-smooth whenever 
$\calz_f$ is. To show that $\calz_g$ is Fano we first show that the adjunction formula 
holds and then compute its index. Now the hypersurface $\calz_g$ has degree $dp$ in 
the weighted projective space $\bbp(d,p\bfw)$ which is not well-formed. The map 
$z_0^p\mapsto z_0$ gives an isomorphism of embeddings of algebraic varieties
$$\matrix{\calz_g&\hookrightarrow &\bbp(d,p\bfw)\cr
                \decdnar{\approx}&&\decdnar{\approx}\cr
                \bbp(\bfw)&\hookrightarrow &\bbp(d,\bfw).\cr}\leqno{\chfol.13}$$
Here the bottom inclusion is the hyperplane $z_0=0,$ and the isomorphism 
$\calz_g\approx \bbp(\bfw)$ is given by the map $z_0^p\mapsto z_0\mapsto z_0-f.$ But 
$\bbp(\bfw)$ is isomorphic as an algebraic variety to $\bbp(\bar{\bfw}).$ This gives an 
isomorphism of algebraic varieties $\calz_g\approx \bbp(\bar{\bfw}).$ Since the Fano index 
is an algebric invariant, we have $I(\calz_g)=|\bar{\bfw}|.$ \hfill\za

\noindent{\sc Remark} \chfol.14: In Section 7 we shall see that the link associated to 
the hypersurface $\calz_g$ of Lemma \chfol.12 is always a rational homology sphere.

%\vfill\eject
\bigskip
\baselineskip = 10 truept
\centerline{\bf \exfol. The K\"ahler Orbifolds Associated to Homotopy Spheres 
}  
\bigskip

In this section we discuss the K\"ahler orbifolds associated to the 
characteristic foliation of the two types of homotopy spheres $\grS_k^{4m-1}$ 
and $\grS_p^{4m+1}$ discussed in section 1. We shall denote by $\grS^{2n+1}$ 
these when there is no need to distinguish between them. We treat the case of 
$\grS_k^{4m-1}$ first.

\noindent{\sc Proposition} \exfol.1: \tensl For all $k\geq 1$ and $m\geq 2,$
the Sasakian structure on  $\grS_k^{4m-1}$ given by Theorem \chfol.5 is 
positive, quasi-regular of order $6(6k-1)$ and the space of leaves $\calz_k$ 
is a Fano orbifold given as a hypersurface in the weighted projective space
$\bbp(6,2(6k-1),3(6k-1),\cdots,3(6k-1))$ cut out by the equation 
$$z_0^{6k-1}+z_1^3+z_2^2+\cdots +z_{2m}^2=0.$$ 
Moreover, the natural projection $\pi_k:\grS^{4m-1}_k\ra{1.3} \calz_k$ is a
principal $S^1$ V-bundle over the Fano orbifold $\calz_k$ with orbifold Fano index 
$I=2m+1.$ The orbifold first Chern class of this V-bundle is  
$(c_1(\calz_k)/I)\in  H^2_{orb}(\calz_k,\bbz)$ which is associated to the 
generator $\calo_{\calz_k}(3(6k-1))$ in $\hbox{Pic}^{orb}(\calz_k).$ \tenrm

\noindent{\sc Proof}: With the exception of the statements
about the order and the Fano condition this follows immediately from the
commutative diagram
$$\matrix{\grS^{4m-1}_k &\ra{2.5}& S^{4m+1}_\bfw&\cr
  \decdnar{\pi}&&\decdnar{} &\cr
   \calz_k &\ra{2.5} &\bbp(\bfw)&\cr}\leqno{\exfol.2}$$
where the weight vector $\bfw$ is given by \exsas.5, the horizontal arrows are
Sasakian and K\"ahlerian embeddings, respectively, and the vertical arrows are
the natural projections. Furthermore, assuming that $\calz_k$ is Fano, it follows by the 
inversion theorem of [BG1] $\grS^{4m-1}_k$ is the total space of the principal $S^1$ 
V-bundle over the orbifold $\calz_k$ whose first Chern class is $c_1(\calz_k)/I.$ So 
we need to compute the index $I.$ First we note that adjunction holds since by the 
computation of the order in Lemma \exfol.2 below, $\calz_k$ is a Cartier divisor. Now by 
Lemma \chfol.12 $\calz_k$ is isomorphic as an algebraic variety to a well-formed weighted 
projective space. This isomorphism is determined by the sequence of transformations
$$(z_0^{6k-1},z_1^3,z_2^2,\cdots,z_{2m}^2)\mapsto 
(z_0,z_1^3,z_2^2,\cdots,z_{2m}^2)\mapsto (z_0,z_1,z_2^2,\cdots,z_{2m}^2).$$
giving the sequence of isomorphisms
$$\bbp(6,2(6k-1),3(6k-1),\cdots,3(6k-1))\approx \bbp(6,2,3,\cdots,3)\approx 
\bbp(2,2,1,\cdots,1).$$ 
Thus, as algebraic varieties the $\calz_k$ are isomorphic for each $k=1,\cdots$ to the 
weighted quadric 
$$z_0+z_1+z_2^2+\cdots +z_{2m}^2=0$$
in $\bbp(2,2,1,\cdots,1)$ which can be transformed to $z_0=0$ giving the 
weighted projective space $\bbp(2,1,\cdots,1).$ There are $2m-1$ ones, so the Fano 
index is $2m+1.$ Furthermore, $a_\bfw=\lcm(d_0,\cdots,d_{2m})=3(6k-1),$ so 
$\calo_{\calz_k}(3(6k-1))$ is the positive generator in $\hbox{Pic}^{orb}(\calz_k).$

The local uniformizing groups can be computed by
analyzing the singular stratum of $\calz_k$ and this is done in the lemma
below where it is shown that the order of $\calz_k$ is $\upsilon(\calz_k)=6(6k-1).$  
The positivity of the Sasakian structure is implied by the Fano 
condition, namely $I>0.$ 
\hfill\za

The computation of the order is given in:

\noindent{\sc Lemma} \exfol.2: \tensl The (orbifold) singular stratum of
$\calz_k$ is the union $Y_0\cup Y_1\cup Y_2$ where 
$$Y_0=\{[\bfz]\in \calz_k|~z_0=0\}, \quad Y_1=\{[\bfz]\in \calz_k|~z_1=0\},
\quad  Y_2=\{[\bfz]\in \calz_k|~z_2=\cdots =z_{2m}=0\}.$$
The local uniformizing groups $\grG_i$ on $Y_i$ are 
$$\grG_0\approx \bbz_{6k-1}, \qquad \grG_1\approx \bbz_3, \qquad \grG_2\approx
\bbz_2,$$  
so the order $\upsilon(\calz_k)$ is $6(6k-1).$ \tenrm

\noindent{\sc Proof}: The condition for determining the fixed point set of
$\grS^{4m-1}_k$ under the action of finite subgroups of the weighted circle
$S^1_\bfw$ generated by the Reeb vector field $\xi_\bfw$ is determined by the
equations
$$\grz^6z_0= z_0,\qquad \grz^{2(6k-1)}z_1=z_1, \qquad
\grz^{3(6k-1)}z_i = z_i~~\hbox{for}~~ i=2,\cdots,2m.\leqno{\exfol.3}$$ 
If $z_0z_1z_i\neq 0$ for some $i=2,\cdots,2m$ there are no fixed points, so we
look at the three cases:
\item{1.} $z_0=0:$ This gives the hypersurface $z_1^3+z_2^2+\cdots
+z_{2m}^2=0$ with isotropy group $\grG_0\approx \bbz_{6k-1}.$
\item{2.} $z_1=0:$ This is another hypersurface $z_0^{6k-1}+z_2^2+\cdots
+z_{2m}^2=0$  with isotropy group $\grG_1\approx \bbz_3.$
\item{3.} $z_i=0$ for all $i=2,\cdots,2m:$ This gives the set of points
determined by $z_0^{6k-1}+z_1^3=0$ with isotropy group $\grG_2\approx \bbz_2.$
\hfill\za

For the case of $\grS_p^{4m+1}$ we have

\noindent{\sc Proposition} \exfol.3: \tensl For all $p\geq 0$ and $m\geq 1,$
the Sasakian structure on  $\grS_p^{4m+1}$ given by Theorem \chfol.5 is 
positive, quasi-regular of order $2p+1$ and the space of leaves $\calz_p$ is a 
Fano orbifold given as a hypersurface in the weighted projective space
$\bbp(2,2p+1\cdots,2p+1)$ cut out by the equation 
$$z_0^{2p+1}+z_1^2+z_2^2+\cdots +z_{2m+1}^2=0.$$ 
Moreover, the natural projection $\pi_k:\grS^{4m+1}_k\ra{1.3} \calz_p$ is a
principal $S^1$ V-bundle over the Fano orbifold $\calz_p$ of orbifold Fano index 
$I=2m+1.$ The orbifold first Chern class of this V-bundle is  
$(c_1(\calz_p)/I)\in  H^2_{orb}(\calz_p,\bbz)$ which is associated to the 
locally V-free sheaf $\calo_{\calz_p}(2p+1).$  \tenrm

\noindent{\sc Proof}: The proof of this proposition together with the following lemma are 
quite similar to that of Proposition \exfol.1 and Lemma \exfol.2. So we only give the 
computation of the index $I.$ In this case as given in Lemma \exfol.4 below, the order 
$\upsilon$ is $2p+1$ and the degree $d=2(2p+1),$ so the order divides the degree, 
and adjunction holds. Here we have $d_0=2p+1,$ and $d_j=1$ for $j>0,$ so 
$a_\bfw=2p+1.$ Thus, $\calo_{\calz_p}(2p+1)$ is the positive generator of 
$\hbox{Pic}^{orb}(\calz_p).$ The isomorphism $\calz_p\approx \bbp(\bar{\bfw})$ is 
given by Lemma \exfol.12 where $\bar{\bfw}=(2,1,\cdots,1)$ with $2m+1$ 
ones, implying that the Fano index of $\calz_p$ is $I=2m+1.$ \hfill\za 

The order $\upsilon(\calz_p)$ is determined by:

\noindent{\sc Lemma} \exfol.4: \tensl The (orbifold) singular stratum of
$\calz_p$ is the quadric $Q_{2m-1}$ given by  
$$z_1^2+\cdots +z_{2m+1}^2=0$$
with the isotropy group $\bbz_{2p+1}.$ So the order $\upsilon(\calz_p)$ is 
$2p+1.$ 
\tenrm

\bigskip
\baselineskip = 10 truept
\centerline{\bf \exproj. Some Exotic Projective Spaces}  
\bigskip

L\'opez de Medrano [LdM] has shown that the quotient space of any fixed point free 
involution on a homotopy sphere is homotopy equivalent to a real projective space 
$\bbr\bbp^n.$ Here we are interested in fixed point free involutions on $\grS^{2n+1}$ 
that  also preserve the Sasakian structure. In particular, we consider the well known 
involution $T$ on $\grS^{4m+1}_p$ defined by $(z_0,z_1,\cdots,z_{2m+1})\mapsto 
(z_0,-z_1,\cdots,-z_{2m+1}).$ One easily sees that this defines a free action of 
$\bbz_2$ on $\grS^{4m+1},$ so by L\'opez de Medrano's theorem the quotient manifold 
$\grS^{4m+1}_p/T$ is homotopy equivalent to $\bbr\bbp^{4m+1}.$ Furthermore, it is clear 
that this $T\approx\bbz_2$ 
is a subgroup of the weighted circle group $S^1_\bfw$ with weights $\bfw 
=(2,2p+1,\cdots,2p+1).$ Thus, the deformation class of Sasakian structures passes to 
the quotient $\grS^{4m+1}_p/T$ to give a deformation class of Sasakian structures 
together with their characteristic foliation $\calf_\xi.$ Furthermore, by Proposition 
\exfol.3 the Sasakian structures are positive, and so by [BGN2] the homotopy real 
projective spaces $\grS^{4m+1}_p/T$ admit Sasakian metrics of 
positive Ricci curvature.  We have arrived at:

\noindent{\sc Proposition} \exproj.1: \tensl For each $(p,m)\in \bbz^+\times \bbz^+$ the 
homotopy projective space $\grS_p^{4m+1}/\bbz_2$ admits a deformation class of 
positive Sasakian structures; hence, $\grS_p^{4m+1}/T$ admits Sasakian 
structures with positive Ricci curvature. \tenrm

It is known [AB,Bro,Gi1-2] that there are at least $2^{2m}$ diffeomorphism types on a 
homotopy real projective space of dimension $4m+1.$ Atiyah and Bott [AB] and Browder 
[Bro] obtained the bound $2^{2m-1}$ which was then extended to $2^{2m}$ by Giffen 
[Gi1-2]. Furthermore, for $0\leq p<p'\leq 2^{2m}.$ the homotopy real projective spaces 
$\grS^{4m+1}_{p'}/T$ and $\grS^{4m+1}_p/T$ are not diffeomorphic.  Thus, 
each of the diffeomorphism types is represented by deformation classes of positive 
Sasakian structures. This proves Theorem C. \hfil\za

\bigskip
\baselineskip = 10 truept
\centerline{\bf \excon. Exotic Contact Structures}  
\bigskip

With the exception of the last two sentences of Theorem C of the Introduction, 
the inequivalences of the contact structures refer to inequivalences of the 
underlying almost contact structures. The material discussed here is well-known 
and taken from Sato [Sa] and Morita [Mo]. An {\it almost contact 
structure} on an orientable manifold $M^{2n-1}$ can be defined as a reduction 
of the oriented orthonormal frame bundle with group $SO(2n-1)$ to the group 
$U(n-1)\times 1.$ Furthermore, the set $A(M)$ of homotopy classes of almost 
contact structures on $M$ is in one-to-one correspondence with the set of 
homotopy classes of almost complex structures on $M\times \bbr$ [Sa], and
when $M=\grS^{2n-1}$ is a homotopy sphere, the latter is known [Mo]
$$A(\grS^{2n-1})=\pi_{2n-1}(SO(2n)/U(n))\cong 
                    \cases{\bbz \oplus \bbz_2 &if $n\equiv 0(\mod 4)$;\cr
                                 \bbz_{(n-1)!} &if $n\equiv 1(\mod 4)$;\cr
                                 \bbz &if $n\equiv 2(\mod 4)$;\cr
                             \bbz_{{(n-1)!\over 2}} &if $n\equiv 3(\mod        
                            4)$;\cr}\leqno{\excon.1}$$                         
                                                 
\noindent{\sc Proof of Theorem B}: The proof for the homotopy spheres $\grS^{4m-1}$ and 
$\grS^{4m+1}$ are somewhat different. The case $\grS^{4m-1}$ uses the proof of 
Theorem 4.1(i) of Morita [Mo]. For each homotopy sphere $\grS^{4m-1}\in bP_{4m}$ the 
underlying 
almost contact structures are distinguished by Morita's invariant $\partial,$ 
and he finds that 
$$\partial(\grS^{4m-1}_k) = \bigl({2\over S_m}-6\bigr)k,\leqno{\excon.2}$$
where $S_m$ is a certain rational number. Moreover, there are an infinite 
number of $k$ that correspond to the same homotopy sphere $\grS^{4m-1}\in 
bP_{4m}.$ Thus, we will have a countably infinite number of distinct 
underlying almost contact structures on the same homotopy sphere as long as 
$S_m\neq {1\over 3}.$ If $S_m={1\over 3}$ we simply use another description of 
elements in $bP_{4m}$ by Brieskorn manifolds, for example that of \exsph.4. This 
completes the proof for this case.

Representing the homotopy spheres $\grS^{4m+1}$ by 
polynomials of the form \exsph.3, we know that $\grS^{4m+1}$ is the standard 
sphere if $p\equiv \pm 1(\mod 8),$ and the Kervaire sphere if $p\equiv \pm 
3(\mod 8)$ which is exotic when $4m+2\neq 2^i-2.$ In this case the Morita invariant is 
simply given in terms of the Milnor number $\mu$ of the link, viz.
$$\partial(\grS^{4m+1}_p)={1\over 2}\mu(\grS^{4m+1}_p)={1\over 2}(p-1). 
\leqno{\excon.3}$$
Thus, as $p$ varies through  $\pm 1(\mod 8)$ or $\pm 3(\mod 8),$ the invariant
$\partial(\grS^{4m+1}_p)$ varies through $0(\mod 4)$ and $3(\mod 4),$ or 
$1(\mod 4)$ and $2(\mod 4),$ respectively. Since for $m>1$ $a(m)\equiv 0(\mod 
4)$ we can choose $p$ such that $\partial(\grS^{4m+1}_p)$ is any number $\mod 
a(m).$ Thus, there are $a(m)$ distinct underlying almost contact structures which is finite 
in this case. To obtain an infinite number of inequivalent contact structures we use the 
recent results of Ustilovsky [Us] who showed by explicitly computing the contact homology 
of Eliashberg, Giventhal, and Hofer (cf. [Eli]) that for the standard sphere  $S^{4m+1}$ 
distinct $p$'s satisfying $p\equiv \pm 1(\mod 8)$ give inequivalent contact structures. 
However, Ustilovsky's proof works equally well for any homotopy sphere $\grS^{4m+1}.$ 
This completes the proof. \hfill\za

\noindent{\sc Remarks} \excon.4: (1) It is interesting to note that for the 
5-sphere $S^5$ ($m=1$ above), Ustilovsky's results distinguish infinitely many 
inequivalent contact structures, whereas, there is a unique underlying almost 
contact structure. (2) It is tempting to adapt the computations in [Us] to the case of the 
homotopy projective spaces $\bbr\bbp^{4m+1}.$ However, in this case complications 
arise, not the least of which is the fact that the simplification of the contact homology 
arising from the fact that all the Conley-Zehnder indices are even no longer holds.

\vfill\eject
\bigskip
\baselineskip = 10 truept
\centerline{\bf \exbr. Homotopy Spheres as Branched Covers}  
\bigskip

It is clear that the Brieskorn spheres described above can be thought of as 
cyclic branched covers of standard spheres. For example, the Kervaire spheres 
$\grS^{4m+1}_p$ given by \exsph.3 is a $p$-fold cyclic cover of $S^{4m+1}$ 
branched over  the Stiefel manifold $V_2(\bbr^{2m+1})$ realized as an $S^1$ 
bundle over an odd complex quadric. In particular, the exotic Kervaire 
sphere $\grS^{4m+1}_3$ can be realized as a 3-fold cover of $S^{4m+1}$ 
branched over the Stiefel manifold $V_2(\bbr^{2m+1})\subset S^{4m+1}.$ 
Similarly, the homotopy spheres $\grS^{4m-1}_k\in bP_{4m}$ are $6k-1$-fold 
covers of $S^{4m-1}$ branched over the Kervaire sphere $\grS^{4m-3}_3.$ In 
both cases the branching locus is a rational homology sphere. More generally 
Savel'ev [Sav] noticed that under certain conditions the cyclic branched cover 
of a rational homology sphere is a homotopy sphere. Here we formulate 
Savel'ev's result in a way that is more convenient for our purpose. We make 
use of the work of Milnor and Orlik [MO]. 

\noindent{\sc Theorem} \exbr.1: \tensl Let $f(z_1,\cdots,z_n)$ be a weighted 
homogeneous polynomial of degree $d$ and weights $\bfw=(w_1,\cdots,w_n)$ in 
$\bbc^n$ with an isolated singularity at the origin. Let $p\in \bbz^+$ and 
consider the link $L_g$ of the equation 
$$g=z_0^p+f(z_1,\cdots,z_n)=0.$$
Write the numbers ${d\over w_i}$ in irreducible form ${u_i\over v_i},$ and 
suppose that $\gcd(p,u_i)=1$ for each $i=1,\cdots,n.$ Then the link $L_g$  is 
a rational homology sphere with weights ${(d,p\bfw)\over \gcd(p,d)}$ and 
degree $\lcm(p,d).$ Moreover, $L_g$ is a homotopy sphere if and only if the 
link $L_f$ is a rational homology sphere. \tenrm 

\noindent{\sc Proof}: Let us briefly recall the construction of the Alexander 
(characteristic) polynomial $\grD_f(t)$ in [MO] associated to a link $L_f$ of 
dimension $2n-3.$ It is the characteristic polynomial of the monodromy map 
$\bbi-h_*:H_{n-1}(F,\bbz)\ra{1.3} H_{n-1}(F,\bbz)$ 
induced by the $S^1_\bfw$ action on the Milnor fibre $F.$ Thus, 
$\grD(t)=\det(t\bbi-h_*).$ Now both $F$ and its closure $\bar{F}$ are homotopy 
equivalent to a bouquet of $n$-spheres, and the boundary of $\bar{F}$ is the 
link $L_f$ which is $n-3$-connected. The Betti numbers  
$b_{n-2}(L_f)=b_{n-1}(L_f)$ equal the number of factors of $(t-1)$ in 
$\grD_f(t).$ Thus, $L_f$ is a rational homology sphere if and only if 
$\grD(1)\neq 0,$ and a homotopy sphere if and only if $\grD_f(1)=\pm 1.$
Following Milnor and Orlik we let $\grL_j$ denote the divisor of $t^j-1.$ Then 
the divisor of $\grD_f(t)$ is given by
$$\hbox{div}~\grD_f~= \prod_{i=1}^n({\grL_{u_i}\over v_i}-1)\leqno{\exbr.2}$$
Using  the relations  $\grL_a\grL_b=\gcd(a,b)\grL_{lcm(a,b)},$ equation 
\exbr.2 takes the form
$$(-1)^n+\sum a_j\grL_j,$$
where $a_j\in \bbz$ and the sum is taken over the set of all least common 
multiples of all combinations of the $u_1,\cdots,u_n.$ The Alexander 
polynomial is then given by
$$\grD_f(t) = (t-1)^{(-1)^n}\prod_j(t^j-1)^{a_j},\leqno{\exbr.3}$$
and 
$$b_{n-2}(L_f)~=~(-1)^n+\sum_j a_j. \leqno{\exbr.4}$$

Now we compute the divisor $\hbox{div}~\grD_g.$ We have
$$\hbox{div}~\grD_g~=~(\grL_p-1)\hbox{div}~\grD_f~=~(\grL_p-1)\bigl((-1)^n+\sum 
a_j\grL_j\bigr)$$$$~=~\sum_j\gcd(p,j)a_j\grL_{\scriptstyle{\lcm(p,j)}} 
-\sum_ja_j\grL_j+(-1)^n\grL_p+(-1)^{n+1}.$$
Since the $j$'s run through all the least common multiples of the set 
$\{u_1,\cdots,u_n\}$ and $\gcd(p,u_i)=1$ for all $i,$ we see that for all $j,$ 
$\gcd(p,j)=1.$ This implies 
$$b_{n-1}(L_g)~=~ \sum_ja_j-\sum_ja_j+(-1)^n+(-1)^{n+1}~=~0.$$
Thus, $L_g$ is a rational homology sphere. Next we compute the Alexander 
polynomial for $L_g.$
$$\grD_g(t) ~=~ (t-1)^{(-1)^{n+1}} (t^p-1)^{(-1)^n} \prod_j(t^{pj}-1)^{a_j}   
(t^j-1)^{-a_j}$$$$=~(t^{p-1}+\cdots +t+1)^{(-1)^n} \prod_j \Bigl({t^{pj-1}+ 
\cdots +t+1\over t^{j-1}+\cdots +t+1}\Bigr)^{a_j}.\leqno{\exbr.5}$$
This gives
$$\grD_g(1)~=~ p^{\grS_ja_j+(-1)^n}~=~1$$
where by \exbr.4 the last equality holds if and only if $L_f$ is a rational 
homology sphere. \hfill\za

An easy way to insure the condition $\gcd(p,u_i)=1$ for all $i=1,\cdots,n$ is 
to assume that $\gcd(p,d)=1$ which we shall do henceforth. 

\noindent{\sc Theorem} \exbr.6: \tensl Let $L_g$ be the link of a branched 
cover of $S^{2n+1}$ branched over a smooth rational homology sphere $L_f$ that is 
the link of an isolated hypersurface singularity given by a weighted homogeneous 
polynomial $f$ of weights $\bfw_f$ and degree $d_f$ as in Theorem \exbr.1. Suppose 
also that $\gcd(p,d_f)=1,$ and that the corresponding hypersurface $\calz_f\subset 
\bbp(\bfw_f)$ is well-formed. Then $L_g$ is a homotopy sphere admitting a Sasakian 
structure $\cals$ with positive Ricci curvature. \tenrm

\noindent{\sc Proof}: Since $\gcd(p,d_f)=1$ the conditions of Theorem \exbr.1 
are satisfied, so $L_g$ is a homotopy sphere. Furthermore, since $\calz_f$ is 
well-formed it follows from Lemma \chfol.12 that the index $I_g=\bfw_f>0,$ 
so the induced Sasakian structure is positive. The result 
then follows by Theorem \exsas.14. \hfill\za

This theorem can be immediately applied to the rational homology spheres 
constructed in [BGN4]. We have

\noindent{\sc Theorem} \exbr.7: \tensl Let $L_g$ and $L_f$ be as above, in 
particular $\gcd(p,d_f)=1,$ and take $L_f$ to be one of the 184 rational 
homology 7-spheres listed in [BGN4]. Then $L_g$ admits a deformation class 
$\gF(p,\bfw)$ of Sasakian structures $\cals\in \gF(p,\bfw)$ containing 
Sasakian metrics with positive Ricci curvature. Furthermore, $L_g$ is one of 
the following:
\item{1.} If $p$ and $d_f$ are both odd, then $L_g$ is diffeomorphic to the 
standard 9-sphere $S^9.$
\item{2.} If $p$ is odd and $d_f$ is even, then for the 10 cases listed in the 
table in [BGN4] with $d_f$ even we have: 
\itemitem{a.} For 
\medskip
\centerline{
\vbox{\tabskip=0pt \offinterlineskip
\def\tablerule{\noalign{\hrule}}
\halign to200pt {\strut#& \vrule#\tabskip=1em plus2em&
     \hfil#& \vrule#& \hfil#& \vrule#\tabskip=0pt\cr\tablerule
&&\omit\hidewidth $\bfw=(w_1,w_2,w_3,w_4,w_5)$\hidewidth&& 
\omit\hidewidth $d_f$\hidewidth&\cr\tablerule
&&$(97,1531,2201,2775,3253)$&&$9856$&\cr\tablerule
&&$(101,439,559,579,619)$&&$2296$&\cr\tablerule
&&$(103,1321,2337,2845,3251)$&&$9856$&\cr\tablerule
&&$(115,341,523,591,727)$&&$2296$&\cr\tablerule}} }
\itemitem{} $L_g$ is diffeomorphic to the standard 9-sphere $S^9$ if $p\equiv 
\pm 1(\mod 8),$ and to the exotic Kervaire 9-sphere $\grS^9$ if $p\equiv \pm 
3(\mod 8).$
\itemitem{b.} For
\medskip
%------------------------------- ss style
\centerline{
\vbox{\tabskip=0pt \offinterlineskip
\def\tablerule{\noalign{\hrule}}
\halign to200pt {\strut#& \vrule#\tabskip=1em plus2em&
     \hfil#& \vrule#& \hfil#& \vrule#\tabskip=0pt\cr\tablerule
&&\omit\hidewidth $\bfw=(w_1,w_2,w_3,w_4,w_5)$\hidewidth&& 
\omit\hidewidth $d_f$\hidewidth&\cr\tablerule
&&$(155,1075,3532,5835,7064)$&&$17660$&\cr\tablerule
&&$(187,2416,8177,10965,19328)$&&$41072$&\cr\tablerule
&&$(221,2416,5491,13617,19328)$&&$41072$&\cr\tablerule
&&$(316,1727,9577,13345,24648)$&&$49612$&\cr\tablerule
&&$(316,2041,6751,15857,24648)$&&$49612$&\cr\tablerule}} }
\itemitem{} $L_g$ is diffeomorphic to the standard 9-sphere $S^9$ if 
$p^2\equiv \pm 1(\mod 8),$ and to the exotic Kervaire 9-sphere $\grS^9$ if 
$p^2\equiv \pm 3(\mod 8).$
\itemitem{c.} For
\medskip
\centerline{
\vbox{\tabskip=0pt \offinterlineskip
\def\tablerule{\noalign{\hrule}}
\halign to200pt {\strut#& \vrule#\tabskip=1em plus2em&
     \hfil#& \vrule#& \hfil#& \vrule#\tabskip=0pt\cr\tablerule
&&\omit\hidewidth $\bfw=(w_1,w_2,w_3,w_4,w_5)$\hidewidth&& 
\omit\hidewidth $d_f$\hidewidth&\cr\tablerule
&&$(49,334,525,668,763)$&&$2338$&\cr\tablerule}} }
\itemitem{} $L_g$ is diffeomorphic to the standard 9-sphere $S^9$ if 
$p^3\equiv \pm 1(\mod 8),$ and to the exotic Kervaire 9-sphere $\grS^9$ if 
$p^3\equiv \pm 3(\mod 8).$
\item{3.} If $p$ is even and $d_f$ is odd, then $L_g$ is diffeomorphic to the 
$$\cases{S^9 &if  $|H_3(L_f,\bbz)|\equiv \pm 1(\mod 8)$;\cr
                 \grS^9 &if $|H_3(L_f,\bbz)|\equiv \pm 3(\mod 8)$.\cr}$$
\noindent However, the 174 examples given in [BGN4] with odd degree all have 
$|H_3(L_f,\bbz)|\equiv 1(\mod 8);$ hence, in this case $L_g$ is always a 
standard sphere.

\noindent In items 1 and 2 above there is a countable infinity of deformation 
classes of Sasakian structures with positive Ricci curvature on each of the
homotopy 9-spheres.  \tenrm 

\noindent{\sc Proof}: The first part follows immediately from Theorem \exbr.6 
by choosing $p$ such that $\gcd(p,d_f)=1.$ To compute the differential 
invariant we use a Theorem of Levine [Lev, Mil2] which asserts that $L_g$ is the 
standard 9-sphere if $\grD_g(-1)\equiv \pm 1(\mod 8)$ and the exotic Kervaire 
sphere if $\grD_g(-1)\equiv \pm 3(\mod 8).$ So we need to compute $\grD_g(-1).$
First, by factoring we can reduce \exbr.5 further:
$$\grD_g(t) ~=~ (t^{p-1}+\cdots +t+1)^{-1} \prod_j \bigl(t^{(p-1)j}+ 
\cdots +t^j+1\bigr)^{a_j},\leqno{\exbr.7}$$
where $j$ ranges through the set of all least common multiples of the $u_i$'s. 

{\it Case} 1: Here we have both $p$ and $d_f$ odd. Since $d_f$ is odd all the 
$u_i$'s are odd. Thus, $j$ in the product in \exbr.7 is always odd. This 
implies that $\grD_g(-1)=1$ so $L_g$ is the standard 9-sphere.

{\it Case} 2: Again $p$ is odd, but now $d$ is even, so both even and odd 
$j$'s occur in \exbr.7. The odd $j$'s contribute a 1 to the product. Thus, we 
have
$$\grD_g(-1)=\prod_{j~\rm{even}}p^{a_j}=p^{\grS_{\rm{even}~j}a_j}.$$
To proceed further we must look at the list of rational homology 7-spheres in 
[BGN4]. There are precisely 10 on the list with even degree. Moreover, there 
are 2 types that occur as indicated by Lemmas 3.4 and 3.12 of [BGN4], and they 
are distinguished by the divisor of their Alexander polynomials. For the first 
type we have $\hbox{div}~\grD =\grL_d-1,$ and 4 of the 10 are of this type, precisely 
the ones occurring in 2a of the Theorem. In this case there is one even $j,$ 
namely $j=d$ and $a_d=1.$ Thus, we have $\grD_g(-1)=p,$ so 2a follows by 
Levine's result. The remaining 6 rational homology 7-spheres with even degree 
are of type two and for these we have
$$\hbox{div}~\grD = 
n(\bfw)\grL_d+\grL_{m_3}-n(\bfw)\grL_{m_2}-1\leqno{\exbr.8}$$
where $n(\bfw)\in \bbz^+,$ and $d=m_2m_3$ with $m_2$ and $m_3$ relatively 
prime. One can see from the definition of $m_2$ and $m_3$ in [BG4] that $m_3$ 
must be even and $m_2$ odd. So if follow from \exbr.8 that 
$\sum_{j~\rm{even}}a_j=n(\bfw)+1.$ Now $n(\bfw)$ can be read off from the 
table in [BGN4], since the order of the third homology group $H_3(L_f,\bbz)$ 
is 
$m_3^{n(\bfw)+1};$ hence, 2b and 2c follow.

{\it Case} 3: Here $p$ is even and $d_f$ is odd. In this case the form of 
\exbr.7 is indeterminant at $t=-1,$ so we need some further manipulation.
We also need to treat the two types of rational homology 7-spheres separately.
For type 1 we have precisely one $j$ which is odd and equal to $d,$ and 
$a_d=1.$ Thus, from \exbr.5 we see that the denominator in the product 
evaluates to 1 at $t=-1$ since $j$ is odd. Then we have
$$\grD_g(-1) ~=~\lim_{t\rightarrow -1} (t^{p-1}+\cdots +t+1)^{-1} (t^{pd_f-1}+ 
\cdots+t+1\bigr) ~=~\lim_{t\rightarrow -1} t^{p(d_f-1)}+\cdots +t^p+1=d_f.$$ 
For rational homology 7-spheres $L_f$ of type 1 the order of $H_3(L_f,\bbz)$ 
is precisely the degree $d_f.$ For rational homology 7-spheres $L_f$ of type 
2, equation \exbr.8 above holds. Moreover, since $d_f$ is odd, so are $m_2$ 
and $m_3.$ So we have 3 odd $j$'s, namely, $d,m_2,m_3$ and 
$a_d=n(\bfw)=-a_{m_2},~a_{m_3}=1.$ Again the terms in the denominator inside 
the product in \exbr.5 evaluate to 1, and we find
$$\eqalign{\grD_g(-1) &=\lim_{t\rightarrow -1}\Bigl({t^{pm_3-1}+\cdots 
+t+1\over t^{p-1}+\cdots +t+1)}\Bigr) \Bigl({t^{pd_f-1}+ \cdots+t+1\over 
t^{pm_2-1}+\cdots +t+1}\Bigr)^{n(\bfw)}\cr &=\lim_{t\rightarrow -1}    
(t^{p(m_3-1)}+\cdots +t^{p}+1) (t^{p(d_f-1)}+\cdots 
+t^{pm_2}+1)^{n(\bfw)}=m_3^{n(\bfw)+1}.}$$ 
But $m_3^{n(\bfw)+1}$ is the order of $H_3(L_f,\bbz)$ for $L_f$ of type 2.  
Combining this with Levine's result mentioned above proves 3 and the theorem. 
\hfill\za
\vfil\eject
\bigskip
\bigskip
\bigskip
\medskip
\centerline{\bf Bibliography}
\medskip
\medskip
\medskip
\medskip
\font\ninesl=cmsl9
\font\bsc=cmcsc10 at 10truept
\parskip=1.5truept
\baselineskip=11truept
\ninerm

\item{[AB]} {\bsc M.F. Atiyah and R. Bott} {\ninesl A Lefschetz fixed point formula for
elliptic complexes. II. Applications}, Ann. of Math. 88 (1968), 451-491.
\item{[Abe1]} {\bsc K. Abe}, {\ninesl Some examples of non-regular almost contact 
structures on exotic spheres}, T\^ohoku Math. J. 28 (1976), 429-435.
\item{[Abe2]} {\bsc K. Abe}, {\ninesl On a generalization of the Hopf fibration, I}, T\^ohoku 
Math. J. 29 (1977), 335-374.
\item{[AE]} {\bsc K. Abe and J. Erbacher}, {\ninesl Nonregular contact structures on 
Brieskorn manifolds}, Bull. Amer. Math. Soc. 81 (1975), 407-409.
\item{[Ba]} {\bsc W. L. Baily}, {\ninesl On the imbedding of V-manifolds in
projective space}, Amer. J. Math. 79 (1957), 403-430.
\item{[BG1]} {\bsc C. P. Boyer and  K. Galicki}, {\ninesl On Sasakian-Einstein
Geometry}, Int. J. Math. 11 (2000), 873-909.
\item{[BG2]} {\bsc C. P. Boyer and  K. Galicki}, {\ninesl 
3-Sasakian manifolds}. {\it Surveys in differential geometry: 
essays on Einstein manifolds}, 123--184, Surv. Differ. Geom.,
VI, Int. Press, Boston, MA, 1999.
\item{[BG3]} {\bsc C. P. Boyer and  K. Galicki}, {\ninesl New Einstein Metrics
in Dimension Five}, J. Diff. Geom. 57 (2001), 443-463; math.DG/0003174.
\item{[BG4]} {\bsc C. P. Boyer and  K. Galicki}, {\ninesl The twistor space of a 3-Sasakian 
manifold}, Int. J. Math. 8 (1997), 31-60.
\item{[BGN1]} {\bsc C. P. Boyer, K. Galicki, and M. Nakamaye}, {\ninesl On the
Geometry of Sasakian-Einstein 5-Manifolds}, submitted
for publication; math.DG/0012041.
\item{[BGN2]} {\bsc C. P. Boyer, K. Galicki, and M. Nakamaye}, {\ninesl On
Positive Sasakian Geometry}, submitted
for publication; math.DG/0104126.
\item{[BGN3]} {\bsc C. P. Boyer, K. Galicki, and M. Nakamaye}, {\ninesl
Sasakian-Einstein Structures on $\scriptstyle{9\#(S^2\times S^3)}$}, to appear in Trans. 
Amer. Math. Soc.; math.DG/0102181.
\item{[BGN4]} {\bsc C. P. Boyer, K. Galicki, and M. Nakamaye}, {\ninesl Einstein Metrics 
on Rational Homology 7-Spheres}, submitted for publication; math.DG/0108113.
\item{[BGP]} {\bsc C. P. Boyer, K. Galicki, and P. Piccinni}, {\ninesl 
3-Sasakian Geometry, Nilpotent Orbits, and Exceptional Quotients}, preprint 
DG/0007184 , to appear in Ann. Global Anal. Geom.
\item{[Bl]} {\bsc D.E. Blair}, {\ninesl Contact Manifolds in Riemannian Geometry}, LNM 509, 
Springer-Verlag, 1976.
\item{[Bla]} {\bsc R. Blache}, {\ninesl Chern classes and Hirzebruch-Riemann-Roch 
theorem for coherent sheaves on complex-projective orbifolds with isolated 
singularities}, 
Math. Z. 222 (1996), 7-27.
\item{[Br]} {\bsc E. Brieskorn}, {Beispiele zur Differentialtopologie von Singularit\"aten}, 
Invent. Math. 2 (1966), 1-14.
\item{[Bro]} {\bsc W. Browder}, {\ninesl Cobordism invariants, the Kervaire invariant 
and 
fixed point free involutions}, Trans. Amer. Math. Soc. 178 (1973), 193-225.
\item{[Ch]} {\bsc J. Cheeger}, {\ninesl Some examples of manifolds of nonnegative 
curvature}, J. Diff. Geom. 8 (1973), 623-628.
\item{[DK]} {\bsc J.-P. Demailly and J. Koll\'ar}, {\ninesl Semi-continuity of
complex singularity exponents and K\"ahler-Einstein metrics on Fano
orbifolds}, preprint AG/9910118, Ann. Scient. Ec. Norm. Sup. Paris 34 (2001), 525-556.
\item{[Dim]} {\bsc A. Dimca}, {\ninesl Singularities and Topology of
Hypersurfaces}, Springer-Verlag, New York, 1992.
\item{[Dol]} {\bsc I. Dolgachev}, {\ninesl Weighted projective varieties}, in
Proceedings, Group Actions and Vector Fields, Vancouver (1981) LNM 956, 34-71.
\item{[Eli]} {\bsc Y. Eliashberg}, {\ninesl Invariants in Contact Topology}, Doc. Math. J. 
DMV, Extra Volume ICM 1998 II, 327-338.
\item{[ElK]} {\bsc A. El Kacimi-Alaoui}, {\ninesl Op\'erateurs transversalement
elliptiques sur un feuilletage riemannien et applications}, Compositio
Mathematica 79 (1990), 57-106.
\item{[Fle]} {\bsc A.R. Iano-Fletcher}, {\ninesl Working with weighted complete
intersections}, Preprint MPI/89-95, revised version in  {\it Explicit
birational geometry of 3-folds},  A. Corti and M. Reid, eds.,
Cambridge Univ. Press, 2000,  pp 101-173.
\item{[Gi1]} {\bsc C.H. Giffen}, {\ninesl Smooth homotopy projective spaces}, Bull. Amer. 
Math. Soc. 75 (1969), 509-513.
\item{[Gi2]} {\bsc C.H. Giffen}, {\ninesl Weakly complex involutions and cobordism of 
projective spaces}, Ann. Math. 90 (1969), 418-432.
\item{[GM]} {\bsc D. Gromoll and W. Meyer}, {\ninesl An exotic sphere with nonnegative
curvature}, Ann. of Math. 100 (1974), 401-406.
\item{[GZ]} {\bsc K. Grove and W. Ziller}, {\ninesl Curvature and symmetry of Milnor
spheres}, Ann. Math. 152 (2002), 331-367.
\item{[Her]} {\bsc H. Hern\'andez-Andrade} {\ninesl A class of compact manifolds with 
positive Ricci curvature}, Differential Geometry, Proc. Sym. Pure Math. 27 (1975), 
73-87.
\item{[Hi]} {\bsc N. Hitchin}, {\ninesl Harmonic spinors}, Advances in Math. 14 (1974), 
1-55.
\item{[Hz]} {\bsc F. Hirzebruch}, {\ninesl Singularities and exotic spheres}, 
S\'eminaire Bourbaki, 1966/67, Exp. 314, Textes des conf\'erences, o.S., 
Paris: Institut Henri Poincar\'e 1967. (reprinted in F. Hirzebruch, Gesammelte 
abhandlungen, band II).
\item{[HZ]} {\bsc F. Hirzebruch and D. Zagier}, {\ninesl The Atiyah-Singer
Theorem and Elementary Number Theory}, Publish or Perish, Inc., Berkeley, 1974.
\item{[JK1]} {\bsc J.M. Johnson and J. Koll\'ar}, {\ninesl K\"ahler-Einstein
metrics on log del Pezzo surfaces in weighted projective 3-space}, Ann. Inst. 
Fourier 51(1) (2001) 69-79.
\item{[JK2]} {\bsc J.M. Johnson and J. Koll\'ar},
{\ninesl Fano hypersurfaces in
weighted projective 4-spaces}, Experimental Math. 10(1) (2001) 151-158.
\item{[Joy]} {\bsc  D. Joyce}, {\ninesl
Compact manifolds with special holonomy},
Oxford Mathematical Monographs, Oxford University Press, Oxford 2000.
\item{[Kaw]} {\bsc Y. Kawamata}, {\ninesl  Abundance theorem for minimal threefolds}, 
Invent. Math. 108, (1992), 229-246.
\item{[KeMi]} {\bsc M.A. Kervaire and J.W. Milnor}, {\ninesl Groups of homotopy 
spheres: I}, 
Ann. of Math. 77 (1962), 504-537.
\item{[KM]} {\bsc J. Koll\'ar, and S. Mori}, {\ninesl Birational Geometry of
Algebraic Varieties}, Cambridge University Press, 1998.
\item{[La]} {\bsc T. Lance}, {\ninesl Differentiable structures on manifolds}, {\it Surveys 
on
Surgery Theory, Vol I}, 73-104, S. Cappell, A. Ranicki, J. Rosenberg, Eds. Princeton 
Univ. 
Press, Princeton N.J. 2000.
\item{[LdM]} {\bsc S. L\'opez de Medrano}, {\ninesl Involutions on Manifolds}, 
Ergebnisse, 
band 59, Springer-Verlag, New York, 1971.
\item{[Lev]} {\bsc J. Levine}, {\ninesl Polynomial invariants of knots of codimension 
two}, Ann. of Math. 84 (1966), 537-554.
\item{[LM]} {\bsc  R. Lutz and C. Meckert}, {\ninesl Structures de contact sur certaines 
sph\`eres exotiques}, C.R. Acad. Sci. Paris S\`er. A-B 282 (1976) A591-A593.
\item{[Mil1]} {\bsc J. Milnor}, {\ninesl On manifolds homeomorphic to the 7-sphere},
Ann. of Math. 64 (1956), 399-405.
\item{[Mil2]} {\bsc J. Milnor}, {\ninesl Singular Points of Complex
Hypersurfaces}, Ann. of Math. Stud. 61, Princeton Univ. Press, 1968.
\item{[MO]} {\bsc J. Milnor and P. Orlik}, {\ninesl Isolated singularities
defined by weighted homogeneous polynomials}, Topology 9 (1970), 385-393.
\item{[Mo]} {\bsc S. Morita}, {\ninesl A Topological Classification of Complex 
Structures on $\scriptstyle{S^1\times \grS^{2n-1}}$}, Topology 14 (1975), 
13-22.
\item{[Na]} {\bsc J. Nash}, {\ninesl Positive Ricci curvature on fibre bundles}, J. Diff. 
Geom. 14 (1979), 241-254.
\item{[Po]} {\bsc W.A. Poor}, {\ninesl Some exotic spheres with positive Ricci curvature}, 
Math. Ann. 216 (1975), 245-252.
\item{[Sat]} {\bsc I. Satake}, {\ninesl The Gauss-Bonnet theorem for
$V$-manifolds}, J. Math. Soc. Japan V.9 No 4. (1957), 464-476.
\item{[Sa]} {\bsc H. Sato}, {\ninesl Remarks Concerning Contact Manifolds}, 
T\^ohoku Math. J. 29 (1977), 577-584.
\item{[Sav]} {\bsc I.V. Savel'ev} {\ninesl Structure of Singularities of a 
Class of Complex Hypersurfaces}, Mat. Zam. 25 (4) (1979) 497-503; English 
translation: Math. Notes 25 (1979), no. 3--4, 258--261.
\item{[SH]} {\bsc S. Sasaki and C.J. Hsu}, {\ninesl On a property of Brieskorn 
manifolds}, 
T\^ohuku Math. J., 28 (1976), 67-78.
\item{[Tak]} {\bsc T. Takahashi}, {\ninesl Deformations of Sasakian structures
and its applications to the Brieskorn manifolds}, T\^ohoku Math. J. 30 (1978),
37-43.
\item{[Us]} {\bsc I. Ustilovsky}, {\ninesl Infinitely Many Contact Structures 
on $\scriptstyle{S^{4m+1}}$}, Int. Math. Res. Notices 14 (1999), 781-791.
\item{[Vai]} {\bsc I. Vaisman}, {\ninesl On the Sasaki-Hsu contact structure of the 
Brieskorn manifolds}, T\^ohoku Math. Journ. 30 (1978), 553-560.
\item{[Wr]} {\bsc D. Wraith}, {\ninesl Exotic spheres with positive Ricci curvature}, J. 
Diff. Geom. 45 (1997), 638-649.
\item{[Y]} {\bsc S. -T. Yau}, {\ninesl Einstein manifolds with zero Ricci
curvature},
Surveys in Differential Geometry VI:
{\it Essays on Einstein Manifolds};
A supplement to the Journal of Differential Geometry, pp.1-14,
(eds. C. LeBrun, M. Wang); International Press, Cambridge (1999).
\item{[YK]} {\bsc K. Yano and M. Kon}, {\ninesl
Structures on manifolds}, Series in Pure Mathematics 3,
World Scientific Pub. Co., Singapore, 1984.
\medskip
\bigskip \line{ Department of Mathematics and Statistics
\hfil January 2002} \line{ University of New Mexico \hfil }
\line{ Albuquerque, NM 87131 \hfil } \line{ email: cboyer@math.unm.edu,
galicki@math.unm.edu, nakamaye@math.unm.edu\hfil} \line{ web pages:
http://www.math.unm.edu/$\tilde{\phantom{o}}$cboyer,
http://www.math.unm.edu/$\tilde{\phantom{o}}$galicki \hfil}
\bye